\documentclass{amsart}
\usepackage[nohug]{diagrams}
\usepackage{amscd}
\usepackage[mathscr]{eucal}
\usepackage{amssymb, amsmath, amsthm}
\usepackage{amsfonts}
\usepackage{latexsym}
\usepackage{tabularx,array}
\usepackage{hyperref}
\usepackage{latexsym}

\newfont{\msam}{msam10}

\newtheorem{theorem}[]{Theorem}
\newtheorem{proposition}[]{Proposition}
\newtheorem{corollary}[]{Corollary}
\newtheorem{lemma}[]{Lemma}

\newtheorem{definition}[]{Definition}

\theoremstyle{definition}
\newtheorem{remark}[]{Remark}

\newtheorem{conjecture}[]{Conjecture}

\let\nc\newcommand

\def\bthm{\begin{theorem}}
\def\ethm{\end{theorem}}
\def\blemma{\begin{lemma}}
\def\elemma{\end{lemma}}
\def\bproof{\begin{proof}}
\def\eproof{\end{proof}}
\def\bprop{\begin{proposition}}
\def\eprop{\end{proposition}}
\def\bcor{\begin{corollary}}
\def\ecor{\end{corollary}}
\def\bconj{\begin{conjecture}}
\def\econj{\end{conjecture}}
\nc{\la}{\label}
\def\Q{\mathbb{Q}}

\def\L {\boldsymbol{L}}

\def\Com{\mathtt{Com}}

\def\DGL{\mathtt{DGLA}}
\def\DGC{\mathtt{DGC}}
\def\cDGC{\mathtt{DGCC}}
\def\DGLC{\mathtt{DGLC}}

\def\DGA{\mathtt{DGA}}

\def\cDGA{\mathtt{DGCA}}
\def\DGPA{\mathtt{DGPA}}

\def\D{{\mathtt D}}
\def\C{\mathcal{C}}

\def\Ho{{\mathtt{Ho}}}
\def\mfa{\mathfrak{a}}
\nc{\Ob}{{\rm Ob}}
\nc{\Hom}{{\rm{Hom}}}
\nc{\Homcont}{{\mathcal{H}om}}
\nc{\HOM}{\underline{\rm{Hom}}}
\nc{\DER}{\underline{\rm{Der}}}
\nc{\END}{\underline{\rm{End}}}
\nc{\bSym}{\mathbf{Sym}}
\nc{\Ext}{{\rm{Ext}}}
\nc{\Rep}{{\rm{Rep}}}
\nc{\DRep}{{\rm{DRep}}}
\nc{\NCRep}{\widetilde{\rm{Rep}}}
\nc{\RAct}{{\rm{RAct}}}
\nc{\bs}{\backslash}
\nc{\ob}{{\tt{Obs}}}
\nc{\CE}{\mathcal{C}}
\nc{\TP}{{T\!P}}
\nc{\nn}{{{\natural} {\natural}}}
\nc{\n}{{{\natural}}}
\nc{\A}{\mathbb A}
\nc{\B}{{\mathrm{B}}}
\nc{\Ba}{\overline{\mathrm{B}}}
\nc{\bC}{\overline{C}}
\nc{\bOmega}{\boldsymbol{\Omega}}
\nc{\bB}{\boldsymbol{B}}
\nc{\EXT}{\underline{\rm{Ext}}}
\nc{\TOR}{\underline{\rm{Tor}}}
\def\H{\mathrm H}

\def\HC{\mathrm{HC}}

\def\rHC{\overline{\mathrm{HC}}}

\nc{\End}{{\rm{End}}}
\nc{\GL}{{\rm{GL}}}
\nc{\gl}{{\mathfrak{gl}}}
\nc{\rgl}{\overline{{\mathfrak{gl}}}}
\nc{\g}{{\mathfrak{g}}}
\nc{\h}{{\mathfrak{h}}}
\nc{\PGL}{{\rm{PGL}}}
\nc{\SL}{{\rm{SL}}}
\nc{\sll}{\mathfrak{sl}}
\nc{\cn}{ \mbox{\rm c\^{o}ne} }
\nc{\PSL}{{\rm{PSL}}}
\nc{\ad}{{\rm{ad}}}
\nc{\Ad}{{\rm{Ad}}}
\nc{\dlim}{\varinjlim}
\nc{\plim}{\varprojlim}
\nc{\colim}{{\tt{colim}}}

\newcommand{\HH}{{\rm{HH}}}

\newcommand{\Sym}{{\rm{Sym}}}

\newcommand{\id}{{\rm{Id}}}

\newcommand{\Tr}{{\rm{Tr}}}

\newcommand{\Ker}{{\rm{Ker}}}

\newcommand{\Coker}{{\rm{Coker}}}

\newcommand{\into}{\,\hookrightarrow\,}

\newcommand{\onto}{\,\twoheadrightarrow\,}

\def\cb{\boldsymbol{\Omega}}
\def\Sm{\mathtt{Sym}}

\def\bs{\backslash}

\newcommand{\rar}{\xrightarrow{}}

\nc{\env}{\mathrm{End}(V)}
\nc{\FT}{\mathcal{C}}

\numberwithin{equation}{section}
\numberwithin{theorem}{section}
\numberwithin{lemma}{section}
\numberwithin{proposition}{section}
\numberwithin{corollary}{section}
\numberwithin{example}{section}
\numberwithin{remark}{section}
\numberwithin{definition}{section}

\def\drm{\mathrm{DR}^{\bullet}}

\nc{\LL}{\mathcal{L}}
\nc{\rH}{\overline{\H}}
\nc{\rHH}{\overline{\HH}}
\nc{\U}{\mathcal U}

\def\arbreBA{\vcenter{\xymatrix@R=2pt@C=2pt{
&&&&\\
&&&*{}\ar@{-}[ul] & \\
&&*{}\ar@{-}[uurr] \ar@{-}[uull] \ar@{-}[d]     &&\\
&&&&
}}}

\def\arbreAB{\vcenter{\xymatrix@R=2pt@C=2pt{
&&&&\\
&*{}\ar@{-}[ur] &&& \\
&&*{}\ar@{-}[uurr] \ar@{-}[uull] \ar@{-}[d]     &&\\
&&&&
}}}

\def\arbreABC{\vcenter{\xymatrix@R=1pt@C=1pt{
&&&&&&\\
&*{}\ar@{-}[ur] &&&&& \\
&&*{}\ar@{-}[uurr] &&&&\\
&&&*{}\ar@{-}[uuurrr] \ar@{-}[uuulll] \ar@{-}[d] &&&\\
&&&&&&
}}}

\def\arbreBAC{\vcenter{\xymatrix@R=1pt@C=1pt{
&&&&&&\\
&&&*{}\ar@{-}[ul] &&& \\
&&*{}\ar@{-}[uurr] &&&&\\
&&&*{}\ar@{-}[uuurrr] \ar@{-}[uuulll] \ar@{-}[d] &&&\\
&&&&&&
}}}

\def\arbreACB{\vcenter{\xymatrix@R=1pt@C=1pt{
&&&&&&\\
&*{}\ar@{-}[ur] &&&&& \\
&&&&*{}\ar@{-}[uull] &&\\
&&&*{}\ar@{-}[uuurrr] \ar@{-}[uuulll] \ar@{-}[d] &&&\\
&&&&&&
}}}

\def\arbreBCA{\vcenter{\xymatrix@R=1pt@C=1pt{
&&&&&&\\
&&&&&*{}\ar@{-}[ul] & \\
&&*{}\ar@{-}[uurr] &&&&\\
&&&*{}\ar@{-}[uuurrr] \ar@{-}[uuulll] \ar@{-}[d] &&&\\
&&&&&&
}}}

\def\arbreCAB{\vcenter{\xymatrix@R=1pt@C=1pt{
&&&&&&\\
&&&*{}\ar@{-}[ur] &&& \\
&&&&*{}\ar@{-}[uull] &&\\
&&&*{}\ar@{-}[uuurrr] \ar@{-}[uuulll] \ar@{-}[d] &&&\\
&&&&&&
}}}

\def\arbreCBA{\vcenter{\xymatrix@R=1pt@C=1pt{
&&&&&&\\
&&&&&*{}\ar@{-}[ul] & \\
&&&&*{}\ar@{-}[uull] &&\\
&&&*{}\ar@{-}[uuurrr] \ar@{-}[uuulll] \ar@{-}[d] &&&\\
&&&&&&
}}}

\def\arbreACA{\vcenter{\xymatrix@R=1pt@C=1pt{
&&&&&&\\
&*{}\ar@{-}[ur] &&&&*{}\ar@{-}[ul] & \\
&&&&&&\\
&&&*{}\ar@{-}[uuurrr] \ar@{-}[uuulll] \ar@{-}[d] &&&\\
&&&&&&
}}}

\begin{document}
\title{Dual Hodge decompositions and derived Poisson brackets}
\author{Yuri Berest}
\address{Department of Mathematics,
Cornell University, Ithaca, NY 14853-4201, USA}
\email{berest@math.cornell.edu}
\author{Ajay C. Ramadoss}
\address{Department of Mathematics,
Indiana University,
Bloomington, IN 47405, USA}
\email{ajcramad@indiana.edu}
\author{Yining Zhang}
\address{Department of Mathematics,
Indiana University,
Bloomington, IN 47405, USA}
\email{yinizhan@indiana.edu}

\begin{abstract}
We study general properties of Hodge-type decompositions of cyclic and Hochschild homology of universal enveloping algebras of (DG) Lie algebras. Our construction generalizes the operadic construction of cyclic homology of Lie algebras due to Getzler and Kapranov \cite{GK}. We give a topological interpretation of such Lie Hodge decompositions in terms of $S^1$-equivariant homology of the free loop space of a simply connected topological space.
We prove that the canonical derived Poisson structure on a universal enveloping algebra arising from a cyclic pairing on the Koszul dual coalgebra preserves the Hodge filtration on cyclic homology. As an application, we show that the Chas-Sullivan Lie algebra
of any simply connected closed manifold carries a natural Hodge filtration. We conjecture that the Chas-Sullivan Lie algebra is actually graded, i.e. the string topology bracket preserves the Hodge decomposition.

\end{abstract}
\maketitle

\section{Introduction and statement of results}
%
It is well known that the cyclic homology of any commutative (DG) algebra $A$ has a natural decomposition
\begin{equation}
\la{HCcomm}
\rHC_\bullet(A) \,\cong\, \bigoplus_{p=1}^{\infty}\, \HC_\bullet^{(p)}(A)\ ,
\end{equation}
which (in the case when $A$ is smooth) arises from a Hodge-style truncation of the de Rham complex of $A$. This decomposition is usually called the Hodge (or  $\lambda$-)decomposition of $ \HC_\bullet(A)$, and it can indeed be linked to classical Hodge theory in algebraic geometry (see  \cite{BV}, \cite{GS}, \cite{L1}  or \cite{L}).

In \cite{bfprw1}, we found that a direct sum decomposition similar to \eqref{HCcomm} exists for the universal enveloping algebra $ \U \mfa $ of any (DG) Lie algebra $ \mfa\,$:
\begin{equation}
\la{hodgeds}
\rHC_{\bullet}(\U\mfa) \,\cong\,\bigoplus_{p=1}^{\infty}\, \HC^{(p)}_{\bullet}(\mfa)\ \text{.}
\end{equation}
The direct summands of \eqref{hodgeds} appeared in \cite{bfprw1} as domains of certain (derived) character maps $\,\Tr_{\mathfrak{g}}(\mfa):\,\HC_\bullet^{(p)}(\mfa) \to \H_\bullet(\mfa, \mathfrak{g})\,$ with values in representation homology of $ \mfa $ in
a finite-dimensional Lie algebra $ \mathfrak{g} $; in terms of these
maps, we gave  a new homological interpretation of the (strong) Macdonald Conjecture for a
reductive Lie algebra $ \mathfrak{g} $ (see  \cite{bfprw1}, Section~9).

The purpose of this paper is two-fold. First, we study general properties of the Hodge-type decomposition \eqref{hodgeds} for an arbitrary Lie algebra $ \mfa $. We extend the construction of \cite{bfprw1} to the Hochschild homology of $\U\mfa $
and show how the resulting homology theories $ \HC_\bullet^{(p)}(\mfa) $ and $ \HH_\bullet^{(p)}(\mfa) $ generalize the operadic construction of cyclic/Hochschild homology for Lie algebras
due to Getzler-Kapranov \cite{GK}.
Furthermore, we give a natural topological interpretation of the  decomposition
\eqref{hodgeds} in terms of Frobenius operations on $S^1$-equivariant homology
$ \H^{S^1}_{\bullet}(\LL X, \Q) $ of the free loop space of a simply connected
topological space $X$.

Our second goal is to clarify the relation of  \eqref{hodgeds}
to a canonical derived Poisson structure on $ \U\mfa $ introduced in \cite{BCER}.
This Poisson structure comes from a cyclic pairing on the Koszul
dual coalgebra of $ \mfa $, and our key observation is that the corresponding Lie
bracket on cyclic homology of $ \U\mfa $ preserves the Hodge filtration associated to
\eqref{hodgeds}. The motivation for studying such cyclic Poisson
structures comes from topology: it is known
that the Chas-Sullivan Lie algebra of a simply connected closed
manifold $M$ is an example of a derived Poisson algebra associated with a natural
cyclic pairing on the Lambrechts-Stanley model of $M$ (see  \cite{BCER}, Section~5.5). In this way, our results reveal a new general property of string topology Lie algebras. We actually expect that the Chas-Sullivan Lie bracket preserves the Hodge
decomposition \eqref{hodgeds} ({\it cf.} Conjecture~\ref{conj1} in Section~\ref{Sect4}).

We now proceed with a detailed discussion of results of the paper. We begin by reviewing the derived functor construction of the Hodge decomposition \eqref{hodgeds} given in \cite{bfprw1}.

Let $k$ be a field of characteristic $0$.  Given a Lie algebra $ \mfa $ over $k$, we consider
the symmetric ad-invariant $k$-multilinear forms on $ \mfa \,$ of a (fixed) degree $ p \ge 1 $. Every such form is induced from the universal
one: $\,\mfa \times \mfa \times \ldots \times \mfa \to \lambda^{(p)}(\mfa) \,$, which takes its values in the space $\,\lambda^{(p)}(\mfa)\,$ of coinvariants of the adjoint representation of $ \mfa $ in $ \Sym^p(\mfa)\,$.
The assignment $\,\mfa \mapsto \lambda^{(p)}(\mfa)\,$ defines a (non-additive) functor on the category of Lie algebras that extends in a canonical way to the category of DG Lie algebras:
\begin{equation}
\la{lam}
\lambda^{(p)}:\,\DGL_k \rar \Com_k \ ,\quad \mfa \mapsto \Sym^p(\mfa)/[\mfa, \Sym^p(\mfa)]\ .
\end{equation}
The category $ \DGL_k $ has a natural model structure (in the sense of Quillen \cite{Q1}), with
weak equivalences being the quasi-isomorphisms of DG Lie algebras. The corresponding homotopy (derived) category $ \Ho(\DGL_k) $ is obtained from $ \DGL_k $ by localizing at the class of
weak equivalences, i.e. by formally
inverting all the quasi-isomorphisms in $ \DGL_k $. The functor \eqref{lam}, however,
does {\it not} preserve quasi-isomorphisms and hence does not descend to the homotopy category
$ \Ho(\DGL_k) $. To remedy this problem, one has to replace $\,\lambda^{(p)}\,$ by its (left) derived functor
\begin{equation}
\la{Llam}
\L\lambda^{(p)}:\,\Ho(\DGL_k) \to \D(k)\ ,
\end{equation}
which takes its values in the derived category $ \D(k)  $
of $k$-complexes. We write  $\,\HC^{(p)}_{\bullet}(\mfa)\,$ for the homology of $\, \L\lambda^{(p)}(\mfa) \,$ and call it the {\it Lie-Hodge homology} of $ \mfa $.

For $ p = 1 $, the functor $ \lambda^{(1)} $ is just abelianization of
Lie algebras; in this case, the existence of $\, \L\lambda^{(1)} \,$ follows from
Quillen's general theory (see \cite[Chapter~II, \S 5]{Q1}), and
$\, \HC^{(1)}_{\bullet}(\mfa) \,$ coincides (up to shift in degree)
with the classical Chevalley-Eilenberg homology $ \H_\bullet(\mfa, k) $ of the Lie algebra $ \mfa $.
For $p=2$, the functor $ \lambda^{(2)} $ was introduced by Drinfeld \cite{Dr}; the existence of $ \L\lambda^{(2)} $ was established by Getzler and Kapranov \cite{GK} who
suggested that $
\HC^{(2)}_{\bullet}(\mfa) $ should be viewed as an (operadic) version of cyclic homology for Lie algebras.
One of the key results of \cite{GK} is the existence of Connes' type periodicity sequence
for Lie cyclic homology, which (in the notation of \cite{GK}, see Section~\ref{GKcomp} below)
reads
\begin{equation}
\la{GKconnes}
\ldots \,\to\, \mathrm{HA}_{n-1}(\mathtt{Lie}, \mfa)\,\to\, {\mathrm{HB}}_{n-1}(\mathtt{Lie}, \mfa) \,\to\, \HC_{n-1}(\mathtt{Lie}, \mfa) \,\to\,\mathrm{HA}_{n-2}(\mathtt{Lie}, \mfa) \,\to\,
\ldots
\end{equation}
Now, for an arbitrary $ p \ge 1 $, the existence of the derived functor $ \L\lambda^{(p)} $ was
proven in \cite{bfprw1} (see {\it op. cit.}, Theorem~7.1), using Brown's Lemma and some  homotopical arguments from \cite{BKR}. 

To construct the direct sum decomposition \eqref{hodgeds}
we observe that each $ \lambda^{(p)} $ comes together with a natural transformation  to the composite functor
$ \U_\n := (\,\mbox{--}\,)_\n \circ \,\U:\,\DGL_k \to \DGA_{k/k} \to \Com_k $, where
$\, (\,\mbox{--}\,)_\n $
denotes the cyclic functor $\,R \mapsto  R/(k+[R,R]) \,$ on the category of (augmented) associative DG algebras. The natural transformations $\,\lambda^{(p)} \to \U_\n \,$ are induced by the symmetrization maps
\begin{equation}
\la{symfun}
\Sym^p(\mfa) \to \U\mfa\ ,\quad x_1 x_2 \ldots x_p\, \mapsto\,
\frac{1}{p!}\,\sum_{\sigma \in {\mathbb S}_p}\, \pm \,x_{\sigma(1)} \cdot x_{\sigma(2)} \cdot
\ldots \cdot x_{\sigma(p)}\ ,
\end{equation}
which, by the Poincar\'e-Birkhoff-Witt Theorem, assemble to an isomorphism of  DG $\mfa$-modules
$\,\Sym_k(\mfa) \cong \U \mfa \,$. From this, it follows that   $\,\lambda^{(p)} \to \U_\n \,$ assemble to an isomorphism of functors
\begin{equation}
\la{eqv1}
\bigoplus_{p=1}^{\infty} \lambda^{(p)} \,\cong \, \U_\n \ .
\end{equation}
On the other hand, by a theorem of Feigin and Tsygan
\cite{FT} (see also \cite{BKR}), the functor $ (\,\mbox{--}\,)_\n $ has a left
derived functor $ \L(\,\mbox{--}\,)_\n:\, \Ho(\DGA_{k/k}) \to \D(k) $ that computes
the reduced cyclic homology $\,\rHC_\bullet(R)\,$ of an associative algebra
$ R \in \DGA_{k/k} $. Since $ \U $ preserves quasi-isomorphisms and maps cofibrant
DG Lie algebras to cofibrant DG associative algebras, the isomorphism \eqref{eqv1}
induces an isomorphism of derived functors from $ \Ho(\DGL_k) $ to $ \D(k) $:
\begin{equation}
\la{eqv2}
\bigoplus_{p=1}^{\infty}\, \L\lambda^{(p)}\, \cong\, \L(\,\mbox{--}\,)_\n \circ \,\U \ .
\end{equation}
At the level of homology, \eqref{eqv2} yields the direct decomposition \eqref{hodgeds} ({\it cf.}~\cite[Theorem 7.2]{bfprw1}).
%
%
As explained in \cite{bfprw1}, the existence of  \eqref{eqv2} is related to the fact that $ \U\mfa $ is a cocommutative Hopf algebra, and in a sense, the Lie Hodge decomposition \eqref{hodgeds} is Koszul
dual to the classical Hodge decomposition \eqref{HCcomm} for commutative algebras.

In the present  paper, we extend the above derived functor construction to the Hochschild homology of $\, \U\mfa \,$:
\begin{equation}
\la{hhoch}
\overline{\HH}_{\bullet}(\U\mfa) \,\cong\,\bigoplus_{p=0}^{\infty}\, \HH^{(p)}_{\bullet}(\mfa)
\ ,
\end{equation}
and relate \eqref{hhoch} to the Hodge decomposition \eqref{hodgeds} of the
cyclic homology of $ \U\mfa $. More precisely, we prove the following result ({\it cf.}
Theorem~\ref{hodgesbi} and Theorem~\ref{connesvsgksbi} in  Section~\ref{S2}).
\bthm
\la{T1}
For any (DG) Lie algebra $\mfa$, the Connes periodicity sequence for $ \U\mfa $
decomposes into a natural direct sum of Hodge components, each of which is itself a long exact sequence:
\begin{equation*}
\ldots \, \to\, \HC_{n-1}^{(p+1)}(\mfa) \,\xrightarrow{B} \,
{\HH}_n^{(p)}(\mfa) \,\xrightarrow{I}\, \HC_n^{(p)}(\mfa) \,\xrightarrow{S} \,
\HC_{n-2}^{(p+1)}(\mfa) \,\to\, \ldots
\end{equation*}
For $p=1$, the above exact sequence is precisely the Getzler-Kapranov exact
sequence \eqref{GKconnes}.
\ethm
Next, we explain a topological meaning of the decompositions \eqref{hodgeds} and \eqref{hhoch}. Our starting point is a fundamental theorem
of Quillen \cite{Q2} that assigns to every $1$-connected topological space
$X$ a DG Lie algebra $ \mfa_X $ over $\Q$, called a Lie model of $X$.  The Lie algebra $ \mfa_X $ determines the rational homotopy type of $X$ and thus determines any homotopy invariant of $X$ defined over $ \Q $. In particular, it is known \cite{M} that the universal enveloping algebra
$ \U \mfa_X $ is quasi-isomorphic to the singular chain algebra
$ C_\bullet(\Omega X, \Q) $ of the based loop space $ \Omega X $
of $X$, while the Hochschild and cyclic homology of $ \U \mfa_X $ are
 isomorphic respectively to the rational homology and rational $S^1$-equivariant homology of the free loop space $ \LL X $ of $X$ (see \cite{Go, J} and also \cite{JM}):
\begin{equation}
\la{jonesiso}
 {\rHH}_{\bullet}(\mathcal U\mathfrak{a}_X) \,\cong\,
 \overline{\H}_{\bullet}(\LL X, \Q)  \ ,\qquad
  \rHC_{\bullet}(\mathcal U\mathfrak{a}_X) \,\cong\,
 {\rH}^{S^1}_{\bullet}(\LL X, \Q) \ .
\end{equation}
In view of \eqref{jonesiso}, it is natural to ask if the Hodge decompositions
\eqref{hodgeds} and \eqref{hhoch} for  $\mfa_X$ can be interpreted in terms of homology of the free loop space $ \LL X$.
To answer this question we recall that  $ \LL X := {\rm Map}(S^1, X) $
carries a natural circle action (induced by the action of  $ S^1 $ on itself).
Considering the $n$-fold covering of the circle:
$\,  S^1 \to S^1 $, $ \, e^{i\theta} \mapsto e^{i n \theta}\,$, we may
define, for each $ n \ge 0 $, $\Q$-linear operations  on the homology
of $\LL X$ (see Section~\ref{Sect4.1} for details):
\begin{equation*}
\Phi_X^n:\,  {\rH}_{\bullet}(\LL X, \Q) \to  {\rH}_{\bullet}(\LL X, \Q)
\ ,\qquad \tilde{\Phi}_X^n :\,  {\rH}^{S^1}_{\bullet}(\LL X, \Q)
\to   {\rH}^{S^1}_{\bullet}(\LL X, \Q)\ .
\end{equation*}
called the Frobenius operations. Now, let ${\rH}^{(p)}_{\bullet}(\LL X, \Q) $ and  $ {\rH}^{S^1,\,(p)}_{\bullet}(\LL X, \Q) $ denote the common eigenspaces of these
Frobenius operations ($ \Phi_X^n $ and $ \tilde{\Phi}_X^n $, respectively), corresponding to the eigenvalues $ n^p $ for all $ n \ge 0 $. Then, we have the following theorem ({\it cf.} Theorem~\ref{top2} in
Section~\ref{Sect4}), which is our second main result in the paper.
\bthm
\la{T2}
For each $p\ge 0$, there are natural  isomorphisms
\begin{equation*}
{\rHH}^{(p)}_{\bullet}(\mathfrak{a}_X) \cong {\rH}^{(p)}_{\bullet}(\mathcal L X, \Q) \ ,\qquad
\rHC^{(p)}_{\bullet}(\mathfrak{a}_X)  \cong
 {\rH}^{S^1, \,(p-1)}_{\bullet}(\LL X, \Q)
\end{equation*}
compatible with the isomorphisms  \eqref{jonesiso}.
\ethm
We prove Theorem~\ref{T2} by dualizing the classical Hodge decomposition \eqref{HCcomm}
for Sullivan's commutative DG algebra model of $ X $ constructed in \cite{BFG}.

In Section~\ref{Sect3}, we turn to derived Poisson structures. This notion was introduced in \cite{BCER} as a natural homological extension of the notion of an $\H_0$-Poisson structure proposed by Crawley-Boevey \cite{CB}.
Roughly speaking, a Poisson structure on an (augmented)
DG algebra $A$ is a DG Lie algebra structure on the cyclic space
$\, A_\n := A/(k+[A,A])\,$ induced by derivations of $A$ (see Section~\ref{Sect3.1}
for a precise definition). We introduce a category of Poisson DG algebras,$\, \DGPA_k \,$,
together with an appropriate class of weak equivalences. Ideally, one would like to make
$ \DGPA_k $ a (closed) model category in order to have a well-behaved homotopy category $ \Ho(\DGPA_k) $ to consider homotopy invariant structures and derived functors on $ \DGPA_k $. Although we do not achieve this goal in the present paper, we observe (see Proposition~\ref{homcat}) that $ \DGPA_k $ has a weaker property of being a (saturated) {\it homotopical category} in the sense of Dwyer-Hirschhorn-Kan-Smith \cite{DHKS}. Thanks to general results of \cite{DHKS}, this property still allows one to define a homotopy category $ \Ho(\DGPA_k) $ that has good formal properties and supports a meaningful theory of derived functors on $ \DGPA_k $. Having defined the homotopy category $ \Ho(\DGPA_k)$ of Poisson DG algebras, we then define a {\it derived Poisson algebra} to be simply an object\footnote{For technical reasons, we will also assume that the associative DG algebra $A$ on which we put a Poisson structure is cofibrant as an object of $ \DGA_{k/k} $.} of $\, \Ho(\DGPA_k)$. The key result here is Proposition~\ref{lieonft} which says that the cyclic
homology of any derived Poisson algebra $A$ carries a well-defined  bracket $\,\{\,\mbox{--}\,,\mbox{--} \,\}: \rHC_\bullet(A) \times \rHC_\bullet(A)
\to \rHC_\bullet(A) \,$, making $ \rHC_\bullet(A) $ a (graded) Lie algebra.

In Section~\ref{Sect3.1}, following \cite{BCER}, we consider a particular class of derived Poisson structures on $ \U\mfa $ that arise from a cyclic pairing on the Koszul dual coalgebra
of $ \mfa $. We call such derived Poisson structures {\it cyclic}. The following theorem clarifies the relation between cyclic Poisson structures and the Hodge decompostion \eqref{hodgeds} of $ \U\mfa $ ({\it cf.} Theorem~\ref{thodgefiltration} in Section~\ref{Sect3}).
\bthm
\la{T3}
The Lie bracket on $ \rHC_\bullet(\U\mfa) $ induced by a cyclic derived Poisson structure
on $ \U\mfa $ preserves the {\rm Hodge filtration:}
\begin{equation*}
F_p\rHC_{\bullet}(\U\mfa)\,:=\,
\bigoplus_{r \leq p+2} \HC^{(r)}_{\bullet}(\mfa) \ ,
\end{equation*}
thus making $ \rHC_\bullet(\U\mfa) $ a filtered Lie algebra. Moreover, in general, we have
$$
\{\HC^{(1)}_{\bullet}(\mfa), \,\HC^{(p)}_{\bullet}(\mfa)\} \,\subseteq \,
\HC^{(p-1)}_{\bullet}(\mfa) \quad \mbox{and}\quad
\{\HC^{(2)}_{\bullet}(\mfa),\, \HC^{(p)}_{\bullet}(\mfa)\} \subseteq \, \HC^{(p)}_{\bullet}(\mfa)\ ,
\quad \forall\,p \geq 1\ .
$$
In particular,  $\HC^{(2)}_{\bullet}(\mfa)$ is a Lie subalgebra of $\rHC_{\bullet}(\mathcal U\mfa)$ and $\rHC_{\bullet}(\mathcal U\mfa)$ is a graded Lie module over $\HC^{(2)}_{\bullet}(\mfa)$.
\ethm
We illustrate Theorem~\ref{T3} in a number of explicit examples, which include Abelian, unimodular and necklace Lie algebras (see Section~\ref{Sect3.5}). In particular, in the case of necklace Lie algebras, we show that
the Lie bracket on $ \rHC_\bullet(\U\mfa) $ induced by a cyclic  Poisson structure  does {\it not} preserve the Hodge decomposition \eqref{hodgeds}. Thus, the result of Theorem~\ref{T3} cannot be strengthened in the obvious way.

Finally, we apply the results of Theorem~\ref{T2} and Theorem~\ref{T3} to the string topology Lie algebra of Chas and Sullivan \cite{ChS}. We recall
that Chas and Sullivan have shown that the $S^1$-equivariant homology
$ {\rH}^{S^1}_{\bullet}(\LL M, \Q) $ of the free loop space of any
smooth compact oriented manifold $M$ carries a natural Lie algebra
structure. Their construction uses  the transversal intersection product of chains which is difficult to realize algebraically
in general ({\it cf.} \cite{CEG}). However, if $M$ is simply connected,
a theorem of Lambrechts and Stanley \cite{LS} provides a
finite-dimensional commutative DG algebra model
for $M$, whose linear dual coalgebra is Koszul
dual to Quillen's Lie model $ \mfa_M $.  The  Lambrechts-Stanley
algebra comes with a natural cyclic pairing which yields a cyclic pairing
on the dual coalgebra. Now, it turns out that the associated cyclic Poisson
structure on $ \U \mfa_M $ induces a Lie bracket on $ \rHC_\bullet(\U\mfa_M) $ that corresponds (under the isomorphism \eqref{jonesiso})
precisely to the Chas-Sullivan bracket on $ {\rH}^{S^1}_{\bullet}(\LL M, \Q) $. As mentioned above, this fact was the main motivation for us to introduce and study derived Poisson structures in general (see \cite{BCER}).

Our last theorem in this Introduction gathers together some properties of the Chas-Sullivan Lie algebras arising from results of the present paper.
({\it cf.} Theorem~\ref{tstringhomology} in Section~\ref{Sect4}).
\bthm
\la{T4}
Let $M$ be a simply connected smooth oriented closed manifold of dimension $d$.

$(i)$
The string topology Lie algebra of  $M$ is filtered as a Lie algebra, with Lie
bracket of degree $ 2-d $ preserving the following {\rm Hodge filtration}
$$
F_p \rH^{S^1}_{\bullet}(\LL M, \Q) \,:=\, \bigoplus_{q \leq p+1} \rH^{S^1,\, (q)}_{\bullet}(\LL M, \Q)\,\text{.}
$$

\noindent
$(ii)$ The  homology $ \rH_{\bullet}(\LL M, \Q)$ of  the free loop space $ \LL M$ is filtered as a Lie module over the string topology Lie algebra of $M$  with respect to the following Hodge filtration
$$ F_p\overline{\H}_{\bullet}(\mathcal L M, \Q) \,:=\, \bigoplus_{q \leq p+2} \overline{\H}^{(q)}_{\bullet}(\mathcal L M, \Q)\,\text{.} $$

\noindent
$(iii)$ The string topology Lie bracket restricts to the first Hodge component $\, {\rH}^{S^1,\,(1)}_{\bullet}(\mathcal L M, \Q)$, making it a Lie algebra. Further, $ {\rH}^{S^1}_{\bullet}(\mathcal L M, \Q)$  is a graded Lie module over $\overline{\H}^{S^1,\,(1)}_{\bullet}(\mathcal L M, \Q)$ with the grading given by the Hodge decomposition of $\,\overline{\H}^{S^1}_{\bullet}(\mathcal L M, \Q)$.
\ethm
We believe that, unlike Theorem~\ref{T3}, Theorem~\ref{T4} may be strengthened. In particular, we expect  that the Chas-Sullivan Lie bracket of a closed  $d$-dimensional manifold preserves not only the Hodge filtration but actually the Hodge decomposition \eqref{hodgeds}, thus making the string topology Lie algebra
a {\it graded} Lie algebra with respect to $p$-degree:
$$
\{\overline{\H}^{S^1,\,(p)}_{\bullet}(\mathcal L M, \Q), \ \overline{\H}^{S^1,\,(q)}_{\bullet}(\mathcal L M, \Q)\}\, \subseteq\, \overline{\H}^{S^1,\, (p+q-1)}_{\bullet}(\mathcal L M, \Q)\ , \quad \forall\,p,q \ge 1\ .
$$
This is part of a more general Conjecture~\ref{conj1} stated in Section~\ref{Sect4}.

\subsection*{Notation}
Throughout this paper, we denote by $\DGA_{k/k}$ (resp., $\cDGA_{k/k}$) the category of augmented, non-negatively graded DG algebras (resp., commutative DG algebras). The  category of non-negatively graded coaugmented, conilpotent DG coalgebras (resp., cocommutative DG coalgebras) will be denoted by $\DGC_{k/k}$ (resp., $\cDGC_{k/k}$). The category of non-negatively graded DG Lie algebras will be denoted by $\DGL_k$.

For $R \in\DGA_{k/k}$, let $R_\n\,:=\, R/(k+[R,R])$.  For an $R$-bimodule $M$, let $M_\n\,:=\, M/[R,M]$. Let $\Omega^1R$ denote the kernel of the multiplication map $R \otimes R \rar R$ (where $R \otimes R$ is equipped with the outer $R$-bimodule structure). The DG $R$-bimodule $\Omega^1R$ represents the complex of derivations $\underline{\mathrm{Der}}(R,M)$. In particular, the universal derivation $\partial\,:\, R \rar \Omega^1R$ is given by $r \mapsto r \otimes 1-1 \otimes r$.

\section{Hodge decomposition for universal enveloping algebras}
\la{S2}
In this section, we construct a Hodge decomposition of  Hochschild homology for the universal enveloping algebras of DG Lie algebras. The corresponding Hodge components are defined in terms of derived functors, similar to the definition of the cyclic Lie-Hodge homology given in~\cite{bfprw1}. The main results of this sections (Theorem~\ref{hodgesbi} and
Theorem~\ref{connesvsgksbi}) explain the relation of our construction to the earlier work of Getzler and Kapranov \cite{GK}.

\subsection{Hodge decomposition of Hochschild homology}
\la{ns1.1}

Let $\mfa$ be a DG Lie algebra, and let $V$ be a right DG $\mfa$-module.
Using the action map on $V$:
$$
V \otimes \mfa \to V\ ,\quad v\otimes x \mapsto v \cdot x\ ,
$$
we define a complex of vector spaces, $\, \theta(\mfa,V) \,$, by factoring  $V \otimes \mfa$ modulo the relations
\begin{equation*}
\la{relunat}
v \otimes [x,y] \,=\, v \cdot x \otimes y- (-1)^{|x||y|} v \cdot y \otimes x \ ,\quad \forall\,v \in V\ ,\ \forall\, x,\,y \in \mfa\ .
\end{equation*}
The action map $\,V \otimes \mfa \to V\,$ then factors through $\theta(\mfa, V)$, giving a canonical
morphism of complexes
\begin{equation}
\la{betav}
\beta_V :\, \theta(\mfa, V) \rar V \,\text{.}
\end{equation}
For example, if $V=k$ is the trivial $\mfa$-module, then
$\,\theta(\mfa, k)\,\cong\,\mfa/[\mfa,\mfa]\,$, and  \eqref{betav} is the zero map. On the other hand,
for $V\,=\, \mathcal U\mfa$ equipped with the (right) adjoint action of $\mfa\,$, the complex $\theta(\mfa, \U\mfa)$ can be identified with $\Omega^1(\U\mfa)_\n\,$, and the map \eqref{betav} is given by
\begin{equation}
\la{betau}
\beta_{\U}:\, \Omega^1(\mathcal U\mfa)_\n \to \U\mfa\ ,\quad
[\alpha \otimes y] \mapsto [\alpha,y]\ ,
\end{equation}
where $\alpha\,\in\,\mathcal U\mfa$, $\,y\,\in\,\mfa$, and $[\alpha \otimes y ]$ denotes the class of  $\alpha \otimes y$ in $\Omega^1(\mathcal U\mfa)_\n$.

Now, fix an integer $\,p\ge1 \,$ and let $ V = \Sym^p(\mfa) $ be the $p$-th symmetric power of $ \mfa $ equipped with the right adjoint action. Write
$ \beta_p(\mfa):\,\theta(\mfa, \Sym^p(\mfa)) \to \Sym^p(\mfa) $ for the corresponding
map \eqref{betav} and define
$$
\phi^{(p)}(\mfa) :=  \mathrm{Cone}[\beta_p(\mfa)]\ .
$$
Note that the functor $\lambda^{(p)}$ introduced in \eqref{lam} is given by
 $ \lambda^{(p)}(\mfa)\,=\,\Coker[\beta_p(\mfa)]\,$; hence, there is a natural transformation of functors
\begin{equation} \la{natI}
 \mathcal I :\, \phi^{(p)} \to \lambda^{(p)} \ ,
\end{equation}
defined by the canonical map $\,\mathrm{Cone}[\beta_p(\mfa)] \to \Coker[\beta_p(\mfa)]\,$.

Next, for $p=0$, we set $\,\phi^{(0)}(\mfa)\,:=\, \theta(\mfa, k)[1]\,$ and define $\phi(\mfa)$ to be the cone of the map \eqref{betau} composed with the natural projection:
$$
\phi(\mfa) := \mathrm{Cone}[\,\Omega^1(\mathcal U\mfa)_\n \xrightarrow{\beta_{\U}} \U\mfa \onto \overline{\U\mfa}\,] \ .
$$
The proof of the following proposition is similar to that of~\cite[Theorem~7.1]{bfprw1}:
we give details for reader's convenience.
\bprop \la{lderived}
The functors $\phi^{(p)}$  and $\phi$ have (total) left derived functors $\L\phi^{(p)}$ and $\,\L\phi\,$ from $\Ho(\DGL_k)$ to $\,\D(k)$.
\eprop
\bproof
First, observe that for any commutative DG algebra $B$, there are natural maps of complexes
$$
\phi^{(p)}(\mfa \otimes B) \rar \phi^{(p)}(\mfa) \otimes B\quad ,\quad 
\phi(\mfa \otimes B) \rar \phi(\mfa) \otimes B\
$$
obtained by extending scalars in the above constructions.
Next, by \cite[Proposition~B.2]{BKR}, if $\mathcal{L}$ is a cofibrant DG Lie algebra, 
two morphisms $f,g\,:\,\mathcal{L} \rar \mfa$ are homotopic in $ \DGL_k $ iff there is a DG Lie algebra 
homomorphism $h\,:\,\mathcal{L} \rar \mfa \otimes \Omega$ such that $h(0)=f$ and $h(1)=g$.
(We call such a homomorphism a {\it homotopy} from $f$ to $g$.). Here $\Omega\,:=\,\Omega(\mathbb A^1_k)$ is the de Rham algebra of the affine line (equipped with homological grading), and $ h(a)$ denotes 
the composite map of $ h $ with $ \id \otimes \mathrm{ev}_a$, where 
$\mathrm{ev}_a\,:\,\Omega \rar k$ is the evaluation map at point $ a \in {\mathbb A}^1_k $.

Now, if $\mathcal{L}$ is a cofibrant DG Lie algebra and $f,g :\,\mathcal{L} \rar \mfa$ are homotopic, with $h$ being a homotopy from $f$ to $g$, one has  the map $ H:\,\phi(\mathcal{L}) \rar \phi(\mfa) \otimes \Omega $ given by the composition
$$
\phi(\mathcal{L})  \xrightarrow{\phi(h)} \phi(\mathfrak{a} \otimes \Omega) \to \phi(\mfa) \otimes \Omega\ .
$$
A straightforward verification shows that $ H(0)=\phi(f)$ and $ H(1)=\phi(g)$. Thus, if $f,g\,:\, \mathcal{L} \rar \mfa$ are homotopic and $\mathcal{L}$ is cofibrant, then $\phi(f)$ and $\phi(g)$ are homotopic as morphisms of complexes.

If $f\,:\,\mathcal{L} \stackrel{\sim}{\to} \mathcal{L}'$ is a weak equivalence between two cofibrant objects 
in $\DGL_k$, by Whitehead's Theorem (see~\cite[Lemma~4.24]{DS}), there exists a map $g\,:\, \mathcal{L}' \rar \mathcal{L}$, such that $fg$ and $gf$ are homotopic to the identities of $ \mathcal{L}' $ and $\mathcal{L}$, respectively. 
It follows that $\phi(fg)$ and $\phi(gf)$ are homotopic to the identities of $\phi(\mathcal{L}')$ and $\phi(\mathcal{L})$, and therefore $\phi(f)$ is a quasi-isomorphism. Hence the functor $\phi$ takes weak equivalences between cofibrant objects to weak equivalences. A similar argument 
shows that the functors $\phi^{(p)}$ take weak equivalences between cofibrant objects to weak equivalences. The existence of  $\,\L\phi\,$ and  $\,\L\phi^{(p)}\,$  follows now from Brown's Lemma in abstract homotopy theory ({\it cf.}~\cite[Lemma 9.9]{DS}).
\eproof
Now, with Proposition~\ref{lderived}, we may define
\begin{equation*}
\la{hhodge}
\HH^{(p)}_{\bullet}(\mfa)\,:=\,
\H_{\bullet}[\L\phi^{(p)}(\mfa)] \ , \quad p\ge 0 \ ,
\end{equation*}
and state first main result of this section:
\bthm \la{hodgehoch}
For any DG Lie algebra $\mfa$, there is a functorial direct sum decomposition
\begin{equation}
\la{hodgeds1}
\overline{\HH}_{\bullet}(\mathcal U\mfa) \,\cong\,\bigoplus_{p=0}^{\infty} \HH^{(p)}_{\bullet}(\mfa)\,\text{.}
\end{equation}
\ethm

By the Poincar\'e-Birkhoff-Witt Theorem (see, e.g., \cite{Q2}, Appendix~B, Theorem~2.3), the  symmetrization maps \eqref{symfun} give a natural isomorphism of right DG $\mfa$-modules $\,
 \bigoplus_{p=0}^{\infty} \Sym^p(\mfa) \cong \U\mfa\,$. This isomorphism
induces an isomorphism of functors $\, \bigoplus_{p=0}^{\infty} \phi^{(p)} \cong \phi\,$,
which, in turn, induces an isomorphism of the corresponding derived functors
from $\Ho(\DGL_k)$ to $\D(k)\,$:
\begin{equation*}
\la{Lphis}
\bigoplus_{p=0}^{\infty} \L\phi^{(p)} \,\cong\,\L\phi\,\text{.}
\end{equation*}
To prove Theorem~\ref{hodgehoch} it thus suffices to prove the following proposition.
\bprop
\la{lrhh}
For any $\mfa\,\in\,\DGL_k$, there is a natural isomorphism $ \H_{\bullet}[\L\phi(\mfa)]\,\cong\, \overline{\HH}_{\bullet}(\mathcal U\mfa)$.
\eprop
\bproof
Let $C\,\in\,\cDGC_{k/k}$ be a cocommutative coalgebra Koszul dual to $\mfa\,\in\,\DGL_k$. 
Assume that $\bar{C}$ is concentrated in strictly positive homological degrees. Then $\,R\,:=\,\cb(C) \stackrel{\sim}{\rar} \U\mfa\,$ is a cofibrant resolution of $ \U\mfa $ in $\DGA_{k/k}$ and $\mathcal L\,:=\, \cb_{\mathtt{Comm}}(C) \stackrel{\sim}{\rar} \mfa $ is a cofibrant resolution
of $\mfa $ in $\DGL_k$. Since $\, R\,\cong\,\,\ U\mathcal L\,$, we have $\theta(\mathcal L, R)\,\cong\,\Omega^1R_\n$. Now, let $ \beta(\mathcal L) $ denote the map \eqref{betau}, with $ \mfa $ replaced by $ \mathcal L $ and composed with the natural projection:
$$ 
\beta(\mathcal L):\, \Omega^1R_\n \rar \bar{R}\,,\,\,\,\,\, \alpha \otimes v \mapsto [\alpha,v]\,\text{.}
$$
Then $\,\L\phi(\mfa) \cong \mathrm{Cone}[\beta(\mathcal L)]\,$ in the derived category
$\D(k)$. Hence, Proposition~\ref{lrhh} is a consequence of the following more general assertion: if $\,R \stackrel{\sim}{\rar} A\,$ be a semi-free resolution of $A$ in $\DGA_{k/k}$, then
$$ 
\overline{\HH}_{\bullet}(A)\,\cong\, \H_{\bullet}(\mathrm{Cone}[\beta(R)])\,\text{.}
$$
To prove this assertion, we will use the approach of \cite[Section 5]{BKR}, which, in turn, is based on Quillen's results \cite{Q}. First, we notice that $ R $ is isomorphic to the tensor algebra $TV$ of a graded vector space $V$; hence, there is an isomorphism of $R$-bimodules  $\,I:\,R \otimes V \otimes R \rar \Omega^1R\,$ given by (see~\cite[Example 3.10]{Q} or~\cite[Section 2.3]{CEEY})
\begin{align*}
& (v_1,\ldots, v_{i-1}) \otimes v_i \otimes (v_{i+1}, \ldots ,v_n) \mapsto  (v_1 ,\ldots v_{i-1},v_i) \otimes (v_{i+1}, \ldots, v_n) - (v_1,\ldots ,v_{i-1}) \otimes  (v_i,v_{i+1} ,\ldots ,v_n) \,\text{.}
\end{align*}
The inverse map $I^{-1}$ induces an isomorphism of graded vector spaces
\begin{align} \la{omeganat}
\Omega^1R_\n\,\cong\, R \otimes V \,\text{.}
\end{align}
 Under \eqref{omeganat}, the map $\bar{\partial}\,:\, \bar{R} \rar \Omega^1R_\n$ becomes\footnote{As explained in~\cite[Appendix A]{BR}, this map may be called the {\it cyclic de Rham differential}. Its kernel is $[\bar{R}, \bar{R}]$ (see {\it loc. cit}).}
\begin{equation}
\la{cyclicderham} \
\bar{\partial}(v_1, \ldots, v_m) \,=\, \sum_{i=1}^m (-1)^{(|v_1|+ \ldots+|v_i| )(|v_{i+1}|+\ldots +|v_m|)} (v_{i+1},\ldots ,v_m,v_1 ,\ldots ,v_{i-1}) \otimes v_i\,\text{.}
\end{equation}
Now, there is a first quadrant bicomplex (see \cite[(5.21)]{BKR})
\begin{equation} 
\la{xc} 
X^{+}(R)\,:=\,
[\, 0 \xleftarrow{}  \bar{R} \xleftarrow{\beta} \Omega^1R_\n \xleftarrow{\bar{\partial}}  \bar{R}   \xleftarrow{\beta}  \ldots\, ] \ , \end{equation}
for which the quotient map
$\,\mathrm{Tot}\,X^+(R) \stackrel{\sim}{\twoheadrightarrow} R_\n\,$
defined by the canonical projection from the first column is a quasi-isomorphism; thus,
$$ 
\H_{\bullet}[\mathrm{Tot}\,X^+(R)]\,\cong\, \rHC_{\bullet}(A)\,\text{.}
$$
By the above identifications, the total complex $\,\mathrm{Tot}\,X_2^+(R) $ of the sub-bicomplex $\,X_2^+(R)\,$ of $\,X^+(R)\,$ comprising of the first two columns is precisely 
$ \mathrm{Cone}[\beta(R)]\,$. On the other hand, by \cite[(5.25)]{BKR}), 
\begin{equation} \la{x2c} \H_{\bullet}[\mathrm{Tot}\,X_2^+(R)]\,\cong\, \overline{\HH}_{\bullet}(A)\,\text{.}\end{equation}
This proves the desired proposition.
\eproof

\subsection{Hodge decomposition of the Connes periodicity sequence}
\la{s1.3}
One of the fundamental properties of cyclic homology is the Connes periodicity exact sequence ({\it cf.}~\cite[2.2.13]{L}).
\begin{equation} \la{connessbi} \begin{diagram}[small] \ldots & \rTo^S & \rHC_{n-1}(A) & \rTo^{B} & \overline{\HH}_n(A) & \rTo^I & \rHC_n(A) & \rTo^S & \rHC_{n-2}(A) & \rTo & \ldots \end{diagram} \,\text{.}\end{equation}
This sequence involves two important operations on cyclic homology: the periodicity operator $S$ and the Connes differential $B$.
\bthm
\la{hodgesbi} Let $\mfa\,\in\,\DGL_k$. The Connes periodicity sequence for $\mathcal U\mfa$ decomposes into a direct sum of Hodge components: the summand of Hodge degree $p$ is given by the long exact sequence
\begin{equation} \la{conneshodgep}  \begin{diagram}[small] \ldots & \rTo^S & \HC_{n-1}^{(p+1)}(\mfa) & \rTo^{B} & {\HH}_n^{(p)}(\mfa) & \rTo^I & \HC_n^{(p)}(\mfa) & \rTo^S & \HC_{n-2}^{(p+1)}(\mfa) & \rTo & \ldots \end{diagram} \,, \end{equation}
with the map $I\,:\,{\HH}_n^{(p)}(\mfa) \rar  \HC_n^{(p)}(\mfa)$ induced on homologies by the natural transformation $\mathcal I\,:\, \L\phi^{(p)}(\mfa) \rar \L\lambda^{(p)}(\mfa)$.
\ethm
In order the prove Theorem~\ref{hodgesbi}, note that $R$ is freely generated by $V\,:=\,\bar{C}[-1]$ as a graded $k$-algebra.  For notational brevity, let $R^{(p)}$ denote the image of $\Sym^p(\mathcal L)$ in $R\,\cong\, \mathcal U\mathcal L$ under the symmetrization map.  Under the isomorphism~\eqref{omeganat}, the direct summand $\theta^{(p)}(\mathcal L)$ of $\Omega^1R_\n$ is identified with $R^{(p)} \otimes V$.
\blemma \la{partialbeta}
For any $p \geq 1$, $\bar{\partial}R^{(p)} \subseteq \theta^{(p-1)}(\mathcal L)$ and $\beta[\theta^{(p)}(\mathcal L)] \subseteq R^{(p)}$.
\elemma
\bproof
The inclusion  $\bar{\partial}R^{(p)} \subseteq \theta^{(p-1)}(\mathcal L)$ follows immediately from Lemma~\ref{cycderhamhodge} proved in the appendix. The inclusion $\beta[\theta^{(p)}(\mathcal L)] \subseteq R^{(p)}$ is a consequence of $V \subseteq \mathcal L$ and $[\Sym^p(\mathcal L), \mathcal L]  \subseteq \Sym^p(\mathcal L)$ in $\mathcal U\mathcal L$ (where we think of the symmetric powers of $\mathcal L$  as subcomplexes of $\mathcal U\mathcal L$ via the symmetrization map).
\eproof
We now proceed with
\begin{proof}[Proof of Theorem~\ref{hodgesbi}]
Let $X^+(R)$ (resp., $X_2^+(R)$) be as in~\eqref{xc} (resp.,~\eqref{x2c}). Recall from~\cite[Section 5]{BKR} that there is an exact sequence of bicomplexes
\begin{equation} \la{xcomplex} 0 \rar X_2^+(R) \rar X^+(R) \rar X^+(R)[2,0] \rar 0\,\text{.} \end{equation}
At the level of total complexes, this gives the exact sequence
\begin{equation} \la{xcomplextot} \begin{diagram}[small] 0 & \rTo &  \mathrm{Tot}\,X_2^+(R) & \rTo^I & \mathrm{Tot}\,X^+(R) & \rTo^S & \mathrm{Tot}\,X^+(R)[2] & \rTo & 0 \end{diagram}\,,\end{equation}
which induces the Connes periodicity sequence on homologies. As an immediate consequence of Lemma~\ref{partialbeta}, we get a direct sum decomposition of bicomplexes
$$X^+(R) \,=\,  \bigoplus_{p=0}^{\infty} X^{+,(p)}(\mfa)\,,$$
where
$$X^{+,(p)}(\mfa)\,:=\,  [ \begin{diagram}[small] 0 & \lTo&  {R}^{(p)} & \lTo^{\beta} & \theta^{(p)}(\mathcal L) & \lTo^{\bar{\partial}} & {R}^{(p+1)}  & \lTo^{\beta} &  \theta^{(p+1)}(\mathcal L) & \lTo^{\bar{\partial}} & {R}^{(p+2)}  & \lTo^{\beta} &  \ldots \end{diagram} ] \,\text{.}$$
In particular, $X_2^+(R) \,=\,  \bigoplus_{p=0}^{\infty} X_2^{+,(p)}(\mfa)\,,$ where $X_2^{+, (p)}(\mfa)$ is the sub-bicomplex of $X^{+,(p)}(\mfa)$ comprising its first two columns. Note that
$$ \mathrm{Tot}\,X_2^{+, (p)}(\mfa) \,\cong\, \L\phi^{(p)}(\mfa) $$
in $\D(k)$, whence $\H_{\bullet}[\mathrm{Tot}\,X_2^{+, (p)}(\mfa)]\,\cong\,\HH^{(p)}_{\bullet}(\mfa)$. Further, since the composite map $\mathrm{Tot}\,X^+(R) \twoheadrightarrow R_\n $
is a quasi-isomorphism, $\mathrm{Tot}\,X^{+,(p)}(\mfa)$ is quasi-isomorphic to $\lambda^{(p)}(\mathcal L)$. It follows that $\H_{\bullet}[\mathrm{Tot}\,X^{+,(p)}(\mfa)]$ is isomorphic to $\HC^{(p)}_{\bullet}(\mfa)$.

In addition, the exact sequence~\eqref{xcomplex} decomposes as a direct sum of Hodge components for $p \geq 0$, with the $p$-th Hodge component given by
\begin{equation} \la{xcomplexhp}  0 \rar X_2^{+,(p)}(\mfa) \rar X^{+,(p)}(\mfa) \rar X^{+,(p+1)}(\mfa)[2,0] \rar 0\,\text{.}  \end{equation}
At the level of total complexes, this gives a Hodge decomposition of the exact sequence~\eqref{xcomplextot}, with the summand in Hodge degree $p$ being
\begin{equation} \la{xcomplexhptot}\begin{diagram}[small] 0 & \rTo &  \mathrm{Tot}\,X_2^{+,(p)}(\mfa) & \rTo^I & \mathrm{Tot}\,X^{+,(p)}(\mfa) & \rTo^S & \mathrm{Tot}\,X^{+,(p+1)}(\mfa)[2] & \rTo & 0 \end{diagram}\,\text{.}\end{equation}
The long exact sequence on homologies corresponding to~\eqref{xcomplexhptot} is~\eqref{conneshodgep}. Finally, since the projection $\mathrm{Tot}\,X^{+,(p)}(\mfa) \twoheadrightarrow \lambda^{(p)}(\mathcal L)$ is a quasi-isomorphism, the natural transformation $\mathcal I\,:\, \L\phi^{(p)}(\mfa) \rar \L\lambda^{(p)}(\mfa)$ is represented by the inclusion of complexes $I$ in~\eqref{xcomplexhptot}. This proves the desired theorem.
\eproof
When $p=0$, the column in degree $0$ of $X_2^{+,(p)}(\mfa)$, namely $R^{(p)}$ vanishes while the column of degree $1$ is isomorphic to $V\,=\,\bar{C}[-1]$. Thus,  $X_2^{+,(0)}(\mfa) \,\cong\, \bar{C}$, whence $\HH^{(0)}_{\bullet}(\mfa)\,\cong\, \H(\mfa;k)$, except in degree $0$ where it vanishes. Since $\HC^{(0)}(\mfa)$ vanishes and since $\HC^{(1)}_{\bullet}(\mfa)\,\cong\,\H_{\bullet+1}(\mfa;k)$, the long exact sequence~\eqref{conneshodgep} for $p=0$ becomes:
$$  \begin{diagram}[small] \ldots &\rTo&\H_n(\mfa;k) & \rTo^{\mathrm{Id}} & {\H}_n(\mfa;k)  & \rTo& 0 & \rTo & \H_{n-1}(\mfa;k) & \rTo^{\mathrm{Id}} & \ldots \end{diagram} \,,$$
with the $B$ map being $\mathrm{Id}$ and the $S$ and $I$ maps vanishing.
\subsection{Comparison to the Getzler-Kapranov Lie cyclic homology}
\la{GKcomp}
We recall that for any cyclic operad $\mathcal P$ and any (DG) $\mathcal P$-algebra $A$, Getzler and Kapranov introduced the (operadic) $\mathcal P$-cyclic homology $\mathrm{HA}_{\bullet}(\mathcal P, A)$ as the homology of the left derived functor of the universal invariant bilinear form on $A$ (see~\cite[Section 4.7, Section 5]{GK}). In addition, they introduced the homologies $\mathrm{HB}_{\bullet}(\mathcal P, A)$ and $\mathrm{HC}_{\bullet}(\mathcal P, A)$ that form the long exact sequence (see~\cite[Section 5.8]{GK})
\begin{equation}
\la{gkseqop}
\ldots \,\to\, \mathrm{HA}_{n}(\mathcal P, A) \,\to\,
{\mathrm{HB}}_{n}(\mathcal P, A) \,\to\, \HC_{n}(\mathcal P, A) \,\to\,
\mathrm{HA}_{n-1}(\mathcal P, A)\,\to\, \ldots
\end{equation}
which we call the Connes periodicity sequence for (operadic) $\mathcal P$-cyclic homology.

It was shown in~\cite[Section 6.10]{GK} that
$$\mathrm{HB}_{n-1}(\mathtt{Lie}, \mfa)\,\cong\, \H_n(\mfa;\mfa)\,,\,\,\,\, \mathrm{HC}_{n-1}(\mathtt{Lie};\mfa)\,\cong\, \H_{n+1}(\mfa;k)\, \text{.}$$
Thus,  $\mathrm{HC}_{n-1}(\mathtt{Lie};\mfa)\,\cong\,\HC^{(1)}_n(\mfa)$. Recall from~\cite[Section 5]{GK} that $\mathrm{HA}_{n-1}(\mathtt{Lie}, \mfa)\,\cong\,\HC^{(2)}_{n-1}(\mfa)$. Our final result in this section is:
\bthm
\la{connesvsgksbi}
For $p=1$, the Hodge component~\eqref{conneshodgep} of the Connes periodicity sequence coincides with the Connes periodicity sequence for (operadic) Lie cyclic homology (see~\cite[Section 6.10]{GK}).
\begin{equation}
\la{gkseq}
\ldots \,\to\, \mathrm{HA}_{n-1}(\mathtt{Lie}, \mfa)\,\to\, {\mathrm{HB}}_{n-1}(\mathtt{Lie}, \mfa) \,\to\, \HC_{n-1}(\mathtt{Lie}, \mfa) \,\to\,\mathrm{HA}_{n-2}(\mathtt{Lie}, \mfa) \,\to\,
\ldots
\end{equation}
\ethm
Before proving Theorem~\ref{connesvsgksbi}, we develop the necessary technical tools.

\subsubsection{}\la{s1.1}

Let  $C\,\in\,\cDGC_{k/k}$ be Koszul dual to $\mfa$. As shown by Quillen~\cite{Q}, $\cb(C)_\n\,\cong\, \overline{\mathrm{CC}}(C)[-1]$, where the right hand side is the reduced cyclic complex of $C$. Thus, $\rHC_{\bullet}(\mathcal U\mfa)\,\cong\, \rHC_{\bullet+1}(C)$. Under this isomorphism, $\HC^{(p)}_{\bullet}(\mfa)$ is identified with the Hodge summand $\rHC_{\bullet+1}^{(p-1)}(C)$ coming from the Hodge decomposition of $\rHC_{\bullet}(C)$ which exists due to cocommutativity of $C$ (see~\cite[Proposition 7.4]{bfprw1}).

In particular, for the rest of this section, let $C$ be the Chevalley-Eilenberg coalgebra $\C(\mfa;k)$. The Hodge decomposition of $\rHC_{\bullet}(C)$ has an explicit description in terms of the de Rham coalgebra of $C$ (see~\cite[Section 2]{bfprw2}). This description is dual to that of the Hodge decomposition of the cyclic homology of a smooth commutative (in fact, symmetric) (DG) algebra $A$ in terms of the de Rham complex of $A$. Explicitly, let $\drm(C)$ denote the mixed de Rham complex of $C$ (see~\cite[Section 2]{bfprw2} for the definition). The negative cyclic complex $\mathrm{CC}^{-}[\drm(C)]$ has a Hodge decomposition due to cocommutativity of $C$. Let $\mathrm{CC}^{-,(p)}[\drm(C)]$ denote the component with Hodge weight $p$. Then (see~\cite{bfprw2}, Proposition~2.2 and Theorem 2.4)
\bprop \la{formderham}
There are natural isomorphisms
$$\rHC_{\bullet}^{(p)}(C)\,\cong\, \H_{\bullet}(\mathrm{CC}^{-,(p)}[\drm(C)])\,\cong\, \H_{\bullet+p}(\ker[d\,:\, \Omega^p_C \rar \Omega^{p-1}_C])\,,$$
where $d\,:\,\Omega^p_C \rar \Omega^{p-1}_C$ is the de Rham differential.
\eprop
It is easy to see that $\Omega^p_C[-p]$ and $\C(\mfa;\Sym^p(\mfa))$ are both isomorphic to  $\Sym^p(\mfa) \otimes \wedge \mfa$ as graded vector spaces. Thus, $\drm(C)$ is isomorphic to $\C(\mfa;\Sym(\mfa))$, which is in turn isomorphic to the (reduced) de Rham algebra $\Omega^{\bullet}_{\Sym(\mfa)}$  as graded vector spaces. We complete this picture with the following
\blemma  \la{chevalley}
Under the above isomorphism, the de Rham differential on $\drm(C)$ is identified with the de Rham differential on $\Omega^{\bullet}_{\Sm(\mfa)/k}$ and the differential  on $\Omega^p_C[-p]$ induced by the differential on $C$ is identified with the Chevalley-Eilenberg differential on $\C(\mfa;\Sym^p(\mfa))$.
\elemma
\bproof
This lemma is formally dual to the following assertion: let $A\,:=\,\C^{\bullet}(\mfa;k)$ be the Chevalley-Eilenberg {\it cochain} complex of $\mfa$ (which is isomorphic as a graded algebra to $\Sym(\mfa^{\ast}[1])$. Then, $\Omega^p_A[-p]\,\cong\,\C^{\bullet}(\mfa;\Sym^p(\mfa^{\ast}))$. For $p=1$, this assertion follows from a direct computation (see~\cite[Section 5.4]{L} for instance). For higher $p$, one notices that the natural map
$\Sym^p_A(\Omega^1_A) \rar \C^{\bullet}(\mfa;\Sym^p(\mfa^{\ast}))$ is compatible with differentials (it induces the $p$-fold cup product) and is an isomorphism of graded $k$-vector spaces.
\eproof
Thus the de Rham differential $d$ on $\Omega_{\Sym(\mfa)/k}$  (anti)commutes with the Chevalley-Eilenberg differential $\delta$ on $\C(\mfa;\Sym(\mfa))$.
 This  makes  $(\Omega^{\bullet}_{\Sym(\mfa)/k}, \delta, d)$ a mixed complex (since $d$ has degree $1$ while $\delta$ has degree $-1$) isomorphic to $\drm(C)$. As a consequence of Proposition~\ref{formderham} and~\cite[Proposition 7.4]{bfprw1},
\bcor \la{liecyclic}
There is a natural isomorphism:
$$ \HC^{(p)}_{\bullet}(\mfa) \,\cong\, \H_{\bullet+1}[\ker(d\,:\,\C(\mfa;\Sym^{p-1}(\mfa)) \rar \C(\mfa;\Sym^{p-2}(\mfa))[-1])] \,\text{.}$$
\ecor
In particular, as shown in~\cite[Section~6.10]{GK},
$$ \HC^{(2)}_n(\mfa)\,\cong\, \H_{n+1}[\ker(d\,:\,\C(\mfa;\mfa) \rar \C(\mfa;k)[-1])]\,\text{.} $$
As a consequence of \eqref{hodgeds} and Corollary~\ref{liecyclic}, we now obtain a different proof of the following result originally due to Kassel~\cite{K} (see also~\cite[Theorem 3.3.7]{L}).
\bcor \la{kassel}
$\rHC_{\bullet}(\mathcal U\mfa)$ is isomorphic to the cyclic homology of the mixed complex $(\Omega^{\bullet}_{\Sym(\mfa)}, \delta, d)$.
\ecor
\bproof
We have natural isomorphisms
\begin{eqnarray*}
\HC_{\bullet}(\Omega^{\bullet}_{\Sym(\mfa)}, \delta, d)\, &\cong &\, \H_{\bullet}[ \Omega^{\bullet}_{\Sym(\mfa)}/d(\Omega^{\bullet}_{\Sym(\mfa)}), \delta]\\
 & \cong & \bigoplus_{p=1}^{\infty} \H_{\bullet}[\mathrm{coker}(d\,:\, \C(\mfa;\Sym^{p+1}(\mfa)) \rar \C(\mfa;\Sym^{p}(\mfa))[-1])]\\
 &\cong &  \bigoplus_{p=1}^{\infty} \H_{\bullet+1}[\ker(d\,:\,\C(\mfa;\Sym^{p-1}(\mfa)) \rar \C(\mfa;\Sym^{p-2}(\mfa))[-1])]\\
 & \cong & \bigoplus_{p=1}^{\infty} \HC^{(p)}_{\bullet}(\mfa)\,\, (\text{by Corollary}~\ref{liecyclic})\\
 & \cong & \rHC_{\bullet}(\mathcal U\mfa)\quad (\text{by}~\eqref{hodgeds})\text{.}\\
\end{eqnarray*}
 where the first isomorphism is due to the fact that $\Omega^{\bullet}_{\Sym(\mfa)}$ is acyclic with respect to the de Rham differential $d$ and the third isomorphism is induced by the de Rham differential.
\eproof

\noindent
\textbf{Remark.} The existence of a Hodge decomposition of $\rHC_{\bullet}(\mathcal U\mfa)$ that is functorial in $\mfa$ follows from Kassel's Corollary~\ref{kassel} and the acyclicity of $\Omega^{\bullet}_{\Sym(\mfa)}$ with respect to the de Rham differential, although this fact seems to have been left unnoticed in~\cite{K, L}. The new interesting fact discovered in~\cite{bfprw1} is the realization of the Hodge components $\HC^{(p)}_{\bullet}(\mfa)$ as homologies of the non-abelian derived functors of universal multilinear forms $\lambda^{(p)}$.

\subsubsection{Proof of Theorem~\ref{connesvsgksbi}} Theorem~\ref{connesvsgksbi} follows from the results of~\cite[Section~6.10]{GK} and the following proposition.
\bprop \la{compkassel}
Let $d^{(p)}\,:\,\C(\mfa;\Sym^{p}(\mfa)) \rar \C(\mfa;\Sym^{p-1}(\mfa))[-1]$ be the de Rham differential. Let $I$ denote the map $d^{(p)}$ thought of as a map from $\C(\mfa;\Sym^{p}(\mfa))$ to $\mathrm{Im}(d^{(p)})$. Then, for all $p \geq 1$, the homology long exact sequence arising from the short exact sequence of complexes
$$ \begin{diagram}[small] 0 & \rTo & \ker(d^{(p)}) & \rTo^B & \C(\mfa;\Sym^{p}(\mfa))  & \rTo^{I} & \ker(d^{(p-1)})[-1] & \rTo & 0 \end{diagram} $$
coincides with the exact sequence~\eqref{conneshodgep}. In particular, we have
$$
\HH^{(p)}_{\bullet}(\mfa)\,\cong\, \H_{\bullet}(\mfa;\Sym^p(\mfa))\ ,
\ \forall\,p \ge 1 \ .
$$
\eprop

\bproof
The bicomplex $X^+(R)$ for $R=\cb(C)$ is formally dual to Tsygan's double complex (in columns of strictly negative degree) for the augmentation ideal $\bar{A}$ of a smooth commutative DGA $A\,\in\,\cDGA_{k/k}$ (with the latter given a bidegree shift of $[-1,1]$. By~\cite[1.4.5 and 2.2.16]{L}, the total complex of this bicomplex is canonically isomorphic to the total complex of Connes' reduced $(b,B)$-bicomplex $\mathcal B^{-}(A)_{\mathrm{red}}[-1,0]$ for the negative cyclic homology $A$. It is easy to see that this isomorphism is compatible with the corresponding $S$-$I$ short exact sequence of complexes for negative cyclic homology. Dually, if $C=\C(\mfa;k)$ (which is a symmetric coalgebra), $X^+(R)$ is canonically isomorphic to Connes' $(b,B)$-complex for $C$ with a shift of $[1,-1]$ (the isomorphism being compatible with the $I$-$S$ short exact sequence of total complexes. By the dual to the Hochschild-Kostant-Rosenberg theorem  and the fact that the de Rham differential on $\Omega^{\bullet}_C$ corresponds to the $B$ differential under the dual HKR map,  there is a dual HKR quasi-isomorphism of bicomplexes
$$ \mathcal B_{\mathrm{DR}}(C)[1,-1] \rar \ \mathcal B(C)_{\mathrm{red}}[1,-1]\,,$$
where
$$
\mathcal B_{\mathrm{DR}}(C)\,:=\,
\left[\ \oplus_p\, \Omega^p_C[-p] \,\xleftarrow{d}\, \oplus_p \,\Omega^{p+1}_C[-p]\, \xleftarrow{d} \,
\ldots \,\right]\,\text{.}
$$
Note that there is a natural Hodge decomposition $\, \mathcal B_{\mathrm{DR}}(C)\,\cong \,\bigoplus_p \mathcal B_{\mathrm{DR}}^{(p)}(C)\,$, where
$$
\mathcal B_{\mathrm{DR}}^{(p)}(C)\,:=\, \left[\,
\Omega^p_C[-p] \,\xleftarrow{d}\,  \Omega^{p+1}_C[-p] \,\xleftarrow{d}\,
\ldots \,\right]\,\text{.} $$
Besides being compatible with the $I$-$S$ short exact sequence of total complexes, the (dual) HKR map is also compatible with Hodge decomposition. The exact sequence~\eqref{conneshodgep} is thus induced by the short exact sequence of complexes
$$
0 \to  \Omega^p_C[-p] \,\xrightarrow{I}\, \mathrm{Tot}\, \mathcal B_{\mathrm{DR}}^{(p)}(C)
\,\xrightarrow{S}\, \mathrm{Tot} \, \mathcal B_{\mathrm{DR}}^{(p+1)}(C)[2] \to 0 \ \text{.}
$$
By Lemma~\ref{chevalley}, $\Omega^p_C[-p]\,\cong\,\C(\mfa;\Sym^p(\mfa))$. Let $\varphi^{(p)}\,:\, \mathrm{Tot}\, \mathcal B_{\mathrm{DR}}^{(p)}(C)  \rar \ker(d^{(p-1)})[-1]$ be the map given by $d^{(p)}$ on $\Omega^p_C[-p]$ and vanishing on all other direct summands. By acyclicity of $\Omega^{\bullet}_C$ with respect to $d$, $\varphi^{(p)}$ is a quasi-isomorphism. Thus, in the derived category $\mathcal D(k)$ of complexes of $k$-vector spaces, there is a commutative diagram where the rows are distinguished triangles and given vertical arrows are isomorphisms.
$$ \begin{diagram}
  \mathrm{Tot} \, \mathcal B_{\mathrm{DR}}^{(p+1)}(C)[1] & \rTo^B & \Omega^p_C[-p] & \rTo^I & \mathrm{Tot}\, \mathcal B_{\mathrm{DR}}^{(p)}(C) & \rTo^S & \mathrm{Tot} \, \mathcal B_{\mathrm{DR}}^{(p+1)}(C)[2] \\
    & & \dTo^{\mathrm{Id}} &  & \dTo^{\varphi^{(p)}} &  &\\
\ker(d^{(p)}) & \rTo^B & \Omega^p_C[-p] & \rTo^I & \ker(d^{(p-1)})[-1] & \rTo & \ker(d^{(p)})[1]\\
\end{diagram}$$
It follows that there is an arrow in $\mathcal D(k)$ completing the above diagram into an isomorphism of distinguished triangles. This proves the desired proposition.
\eproof

\section{Derived Poisson structures}
\la{Sect3}

 In this section, we study derived  Poisson structures on (the universal enveloping algebra of) a DG Lie algebra $\mfa$. We focus on the case when the derived Poisson structure on $\mathcal U\mfa$ arises from a cyclic pairing on a cocommutative DG coalgebra $C$ that is Koszul dual to $\mfa$. Our main observation is that this derived Poisson structure is compatible with Lie Hodge decomposition of $\mathcal U\mfa$.

\subsection{Derived Poisson algebras}
\la{Sect3.1}
We begin by reviewing the notion of a derived Poisson algebra introduced in ~\cite{BCER}. This notion is a higher homological extension of the notion of an ${\rm H}_0$-Poisson algebra proposed by Crawley-Boevey \cite{CB}.

\subsubsection{Definitions}
\la{Defs}
Let $ A$ be an (augmented) DG algebra. The space  $\DER(A)$  of graded $k$-linear derivations of $A$
is naturally a DG Lie algebra with respect to the commutator bracket. Let $\DER(A)^\n $ denote the subcomplex of  $\DER(A) $
comprising derivations with image in $\,k+[A,A] \subseteq A \,$. It is easy to see that $ \DER(A)^\n $ is a DG Lie ideal of $ \DER(A) $,
so that $\,\DER(A)_{\natural} := \DER(A)/\DER(A)^\n $ is a DG Lie algebra.  The natural action of $ \DER(A) $ on $A$ induces a Lie algebra
action of $ \DER(A)_\n $ on the quotient space $ A_\n := A/(k+[A,A]) $. We write $\, \varrho:  \DER(A)_\n \to \END(A_\n) \,$ for the
corresponding DG Lie algebra homomorphism.

Now,  following \cite{BCER},  we define a {\it Poisson structure} on $A$ to be  a DG Lie algebra structure on $\,A_\n \,$ such that the adjoint representation
$ \mbox{\rm ad}:\, A_\n \to \END(A_\n) $ factors through $ \varrho \,$: i.~e., there is a morphism of DG Lie algebras $\,\alpha :\, A_\n \rar \DER(A)_{\natural}\,$ such that $\,\mbox{\rm ad} = \varrho \circ \alpha \,$.  It is easy to see that if $A$ is a commutative DG algebra, then a Poisson structure on $A$ is
the same thing as a (graded) Poisson bracket on $A$. On the other hand, if $A$ is an ordinary $k$-algebra (viewed as a DG algebra), then a Poisson structure on
$A$ is precisely  an ${\rm H}_0$-Poisson structure in the sense of \cite{CB}.

Let $A$ and $B$ be two Poisson DG algebras, i.e. objects of $\DGA_{k/k}$ equipped with Poisson structures.
A {\it morphism} $\,f:\, A \rar B $ of Poisson algebras is then a morphism $ f: A \to  B $ in $\DGA_{k/k} $ such that $ f_{\natural}:\, A_{\natural} \rar B_{\natural} $ is a morphism of DG Lie algebras. With this notion of morphisms, the Poisson DG algebras form a category which we denote $\mathtt{DGPA}_k $.
Note that $\mathtt{DGPA}_k $ comes  with two natural functors: the forgetful functor $ U:\, \mathtt{DGPA}_k  \to \DGA_{k/k} $ and the cyclic functor
$ (\,\mbox{--}\,)_\n : \mathtt{DGPA}_k  \to \DGL_k $. We say that a morphism  $ f $ is a {\it weak equivalence} in $ \mathtt{DGPA}_k $
if $ Uf $ is a weak equivalence in $ \DGA_{k/k} $ and $ f_\n $ is a weak equivalence in $ \DGL_k $; in other words, a weak equivalence in $\mathtt{DGPA}_k $
is a quasi-isomorphism of DG algebras,   $ f: A \to B \,$,  such that the induced map  $ f_{\n}\,:\,A_\n \rar B_\n $ is a quasi-isomorphism of DG Lie algebras.

Although we do not know at the moment whether the category $ \mathtt{DGPA}_k $ carries a Quillen model structure (with weak equivalences specified above), it has a weaker property of being a homotopical category  in the sense of Dwyer-Hirschhorn-Kan-Smith \cite{DHKS}. This still allows one to define a well-behaved homotopy category of Poisson algebras and consider derived functors on $ \mathtt{DGPA}_k $.

\subsubsection{Homotopical categories}
Recall ({\it cf.} \cite{DHKS}) that a {\it homotopical category} is a category $\mathscr{C}$ equipped with a class of morphisms $\mathscr{W} $ (called weak equivalences) that contains all identities
of $ \mathscr{C} $ and satisfies the following  {\it  2-of-6 property}: for every composable triple of morphisms $ f,g,h $ in $\mathscr{C}$, if $ g  f $ and $    h g $ are in $ \mathscr{W} $, then so are $ f, g, h $ and $h g f $.  The 2-of-6 property formally implies, but is stronger than, the usual 2-of-3 property: for every composable pair of morphisms $ f, g $ in $\mathscr{C}$,
if any two of $f$, $g$ and $ gf $ are in $ \mathcal W $,  so is the third.
The class of weak equivalences thus forms a subcategory which contains all objects and all isomorphisms of $ \mathscr{C} $.
The isomorphisms in any category satisfy the 2-of-6 property: indeed, if $ f,g,h $ is a composable triple such that  $ g  f $ and $ h g $ are isomorphisms, then $ g $ has the right inverse
$ f(gf)^{-1} $, which must also be a left inverse, since $ g $ is monic (because $ h g $ is an isomorphism); hence $ g $ and therefore also $f$, $h$ and $hgf$ are isomorphisms.
Thus, any category can be viewed as a homotopical category by taking the weak equivalences to be the isomorphisms\footnote{In \cite{DHKS}, these are called {\it minimal} homotopical categories.}.
Furthermore, any model category is a homotopical category. This follows from the important property
of model categories (see \cite[Proposition~5.1]{Q1}) that the class $ \mathscr{W} $ of weak equivalences is {\it saturated} in $ \mathscr{C} \,$: i.e, it comprises {\it all} the arrows of $ \mathscr{C}$ that become isomorphisms in the localized category
$ \mathscr{C}[ \mathscr{W}^{-1}] \,$. Since the isomorphisms
in $ \mathscr{C}[ \mathscr{W}^{-1}] \,$ satisfy the 2-of-6 property, it follows immediately
that the weak equivalences in a saturated category satisfy the 2-of-6 property.
The category $ \mathscr{C}[ \mathscr{W}^{-1}] \,$ is called the {\it homotopy category} of $ \mathscr{C} $ and usually denoted $ \Ho(\mathscr{C})$. It is a domain
(and target) of homotopical functors and other homotopical structures that constitute the homotopy theory associated to a model category $ \mathscr{C} $.
Now, a key observation of \cite{DHKS} is that a well-behaved homotopy category, including
a meaningful notion of derived functors, can be defined for any homotopical category in which the class of weak equivalences is saturated. Such categories are called in \cite{DHKS} the saturated homotopical categories.

After these preliminaries, we can state our proposition.
\bprop
\la{homcat}
$\,\mathtt{DGPA}_k$ is a saturated homotopical category.
\eprop
\bproof
It suffices to prove that the class of weak equivalences in $\,\mathtt{DGPA}_k$ is saturated. Let  $ \gamma: \mathtt{DGPA}_k \to  \Ho(\mathtt{DGPA}_k)  $ denote the (formal) localization of $ \mathtt{DGPA}_k $ at the class of weak equivalences (similarly, abusing notation, we will write $ \gamma $ for the localizations of $ \DGA_{k/k} $ and $ \DGL_{k} $.)
Since both the forgetful functor $U: \mathtt{DGPA}_k \to \DGA_{k/k} $ and the cyclic functor
$ (\,\mbox{--}\,)_\n:  \mathtt{DGPA}_k \to \DGL_k \to$ preserve weak equivalences, by the universal property of localization, they factor through $ \gamma $, inducing $ \gamma U:\,\Ho(\mathtt{DGPA}_k) \to  \Ho(\DGA_{k/k}) $ and $\,\gamma  (\,\mbox{--}\,)_\n : \Ho(\mathtt{DGPA}_k) \to  \Ho(\DGL_{k}) \,$. Now, if $ f $ is a morphism in  $ \mathtt{DGPA}_k $
such that $ \gamma f $ is an isomorphism in $\Ho(\mathtt{DGPA}_k)$, then $  \gamma (U f)  $  is an isomorphism in $ \Ho(\DGA_{k/k}) $ and
$ \gamma(f_\n) $ is an isomorphism in $ \Ho(\DGL_{k}) $. Since both
$ \DGA_{k/k} $ and $ \DGL_{k} $ are model categories,
the classes of their weak equivalences are saturated. It follows that
$ U f $ and $ f_\n $ are weak equivalences in $ \DGA_{k/k} $ and $ \DGL_{k} $,
respectively. Then, by definition, $ f $ is a weak equivalence in $ \mathtt{DGPA}_k $, so
that the class of weak equivalences  $ \mathtt{DGPA}_k $ is saturated.
\eproof

\subsubsection{Homotopy category of Poisson algebras}
If $ \mathscr{C} $ is a model category, there are two ways to define
a homotopy category of $ \mathscr{C} $: first, one can simply put
$\, \Ho(\mathscr{C}) :=  \mathscr{C}[\mathscr{W}^{-1}]\,$, or alternatively,
one can consider the full subcategory $ \mathscr{C}^{\rm c} $ of cofibrant/fibrant objects
in $ \mathscr{C} $ and then define $ \Ho(\mathscr{C}^{\rm c}) :=  \mathscr{C}^{\rm c}/\!\!\sim $
to be the quotient category of $ \mathscr{C}^{\rm c} $ modulo an
appropriate homotopy equivalence relation\footnote{The category
$ \Ho(\mathscr{C}^{\rm c}) $ is often called the {\it classical homotopy category}
of $ \mathscr{C} $.}. By (an abstract version of) Whitehead's Theorem, the
two definitions are equivalent, the equivalence
$ \Ho(\mathscr{C}^{\rm c}) \stackrel{\sim}{\to} \Ho(\mathscr{C}) $ being induced
by the natural functor $ \mathscr{C}^{\rm c} \into \mathscr{C} \to \Ho(\mathscr{C}) $.
Each definition has its advantages: the first one is more natural (it gives a characterization
of $ \Ho(\mathscr{C}) $ in terms of a universal property and shows that
$ \Ho(\mathscr{C}) $ depends only on the class of weak
equivalences in $ \mathscr{C} $)); the second is more concrete and accessible for computations
(it implies, in particular, that $ \Ho(\mathscr{C}) $ is a locally small
category if so is $\mathscr{C}$.)

For the homotopical category $\mathtt{DGPA}_k $, we
can also define a homotopy category in two ways. First, we can simply put
$$
\Ho(\mathtt{DGPA}_k) :=  \mathtt{DGPA}_k[\mathscr{W}^{-1}]\ ,
$$
where $ \mathscr{W} $ is the class of weak equivalences specified in Section~\ref{Defs}.
Proposition~\ref{homcat} then ensures that $ \Ho(\mathtt{DGPA}_k) $ has properties
similar to those of the homotopy category of a model category (see~\cite[Sect.~33.8]{DHKS}).

Alternatively, following \cite{BCER}, we can mimick the definition of the classical homotopy category of a model category and define\footnote{This category was denoted
$ \Ho^*(\mathtt{NCPoiss}_k) $  in \cite{BCER} and simply referred to as the homotopy
category of Poisson algebras.}
\begin{equation}
\la{defho}
\Ho(\mathtt{DGPA}^{\rm c}_k) :=  \mathtt{DGPA}^{\rm c}_k/\!\!\sim
\end{equation}
The objects in this quotient category are the homotopy classes of Poisson
algebras $A $ whose underlying DG algebras $ U(A) $ are cofibrant
as objects in $ \DGA_{k/k} $. The equivalence relation $\, \sim \,$
is based on the notion of $P$-homotopy (`polynomial homotopy') for Poisson
algebras introduced in \cite{BCER}. We recall that
two morphisms $\,f,g: A \rar B\,$ in $\mathtt{DGPA}_k $ are called {\it $P$-homotopic}
if there is a morphism $\, h: A \rar B \otimes \Omega \,$ such that $ h(0)=f $ and $ h(1)=g $, where $ \Omega = \Omega({\mathbb A}_k^1) $ is the de Rham algebra of the affine line and
$ B \otimes \Omega $ is given the structure of a Poisson DG algebra via the extension of
scalars from $ B $ ({\it cf.}~\cite{BCER}, Sect.~3.1). It is easy to check that the $P$-homotopy defines an equivalence relation on $\, \Hom_{\mathtt{DGPA}}(A,B)\,$ for any objects $A$ and $B$ in $ \mathtt{DGPA}_k $. Now, as in the case of model categories, we have

\bprop
\la{hoeq}
There is a natural functor $\, \Ho(\mathtt{DGPA}^{\rm c}_k) \to \Ho(\mathtt{DGPA}_k) \,$.
\eprop
\bproof
Recall that, for $ B\,\in\,\mathtt{DGPA}_k$, a Poisson DG algebra structure on
$ B \otimes \Omega $ is given by extension of scalars, using the natural isomorphism
$ (B \otimes \Omega)_\n \cong  B_\n \otimes \Omega\, $ ({\it cf.}~\cite[Section 3.1.1]{BCER}).
Then, the inclusion $i\,:\, B \rar B \otimes \Omega$ is a weak equivalence of Poisson algebras: indeed, $\,Ui$ is a weak equivalence of DG algebras, since the de Rham algebra $\Omega$ is acyclic, and $i_\n$ is a weak equivalence of DG Lie algebras, since $i_\n$ can be identified with the inclusion  $B_\n \hookrightarrow B_\n \otimes \Omega$. Thus, if $f,g\,:\, A \rar B$ are $P$-homotopic with a homotopy $h\,:\,A \rar B \otimes \Omega$, then $f=g=i^{-1} \circ h$ in $\Ho(\mathtt{DGPA}_k)$. This proves the desired proposition.
\eproof
\begin{definition}
By a {\it derived Poisson algebra} we mean a cofibrant associative DG algebra $A$ equipped with a
Poisson structure  (in the sense of Defintion~\ref{Defs}), which is viewed up to weak equivalence,
i.e. as an object in $ \Ho(\mathtt{DGPA}_k) $.
\end{definition}

\vspace{1ex}

The above definition differs from that of \cite{BCER}, where
the derived Poisson algebras were simply defined to be the objects of
$ \Ho(\mathtt{DGPA}^{\rm c}_k) $. Proposition~\ref{hoeq} shows, however, that any
derived Poisson algebra in the sense of~\cite{BCER} gives naturally a derived Poisson
algebra in our current sense. The point is that all results of \cite{BCER} established
with the use of explicit $P$-homotopies can be strengthened and reproved in a more
natural way with a weaker notion of equivalence. For example, we have

\bprop
\la{lieonft}
The cyclic homology $ \rHC_{\bullet}(A) $ of any derived Poisson algebra $A$
carries a natural structure of a graded Lie algebra.
\eprop
\bproof
By definition, the functor $\,(\,\mbox{--}\,)_\n:\,
\mathtt{DGPA}^{\rm c}_k \into \mathtt{DGPA}_k \to \DGL_k \,$ preserves weak equivalences,
and hence induces
\begin{equation}
\la{FTf}
(\,\mbox{--}\,)_\n: \, \Ho(\mathtt{DGPA}_k)^{\rm c} \to \Ho(\DGL_k)\ ,
\end{equation}
where $ \Ho(\mathtt{DGPA}_k)^{\rm c} $ denotes the full subcategory of the homotopy
category $ \Ho(\mathtt{DGPA}_k) $ whose objects are cofibrant DG algebras.
Now, the image of $A$ under \eqref{FTf} is a DG Lie algebra $ A_\n $ whose underlying complex
computes the (reduced) cyclic homology of $A$. This is a consequence of \cite{BKR}, Theorem~3.1,
(which is essentially due to Feigin and Tsygan) that states that the functor
$ (\,\mbox{--}\,)_\n $ has a well-defined derived functor $ \L(\,\mbox{--}\,)_\n $ on the category
of DG algebras whose homology agrees with cyclic homology. On cofibrant DG algebras, the values of
$ (\,\mbox{--}\,)_\n $ and $ \L(\,\mbox{--}\,)_\n $ are naturally isomorphic, hence
$ \H_{\bullet}(A_\n) \cong \rHC_{\bullet}(A) $ for any cofibrant $A$. On the other hand,
since $A_\n$ is a DG Lie algebra, $\, \H_{\bullet}(A_\n) $ carries a graded Lie algebra structure. Identifying $ \H_{\bullet}(A_\n) $ with $ \rHC_{\bullet}(A) $ for a cofibrant $A$,
we get a graded Lie algebra structure on $ \rHC_{\bullet}(A) $ claimed by the proposition.
\eproof

Another important result of \cite{BCER} that holds for the derived Poisson algebras in
$ \Ho(\mathtt{DGPA}_k) $ and that motivates our study of these objects is the following

\begin{theorem}[{\it cf.} \cite{BCER}, Theorem~2]
\la{t3s2int}
If $A$ is a derived Poisson DG algebra, then, for any finite-dimensional vector space $ V $,
there is a unique graded Poisson bracket on the representation homology
$ \H_\bullet(A, V)^{\GL(V)} $ of $A$ in $V$, such that the derived character map
$\,
\Tr_V(A)_\bullet:\, \rHC_\bullet(A) \to \H_\bullet(A,V)^{\GL(V)} $
is a Lie algebra homomorphism.
\end{theorem}
We will not reprove Theorem~\ref{t3s2int} here; instead, in Section~\ref{s4}, we will give a generalization of this result to representation homology of Lie algebras.

\subsection{Cyclic pairings} \la{s2.1}
We now describe our basic construction of derived Poisson structures associated with cyclic coalgebras. Recall ({\it cf.}~\cite{GK}) that a graded associative $k$-algebra is called $n$-{\it cyclic} if it is equipped with a symmetric bilinear pairing $\langle \mbox{--},\mbox{--} \rangle\,:\, A \times A \rar k$ of degree $n$ such that
$$ \langle a, bc \rangle \,=\, \pm \langle ca, b \rangle \,,\,\,\,\,\, \forall\,\, a,b,c\,\in\,A\,$$
the signs being determined by the Koszul sign rule. Dually, a graded coalgebra $C$ is called $n$-{\it cyclic} if it is equipped with a symmetric bilinear pairing $ \langle \mbox{--},\mbox{--} \rangle\,:\,  {C} \times {C} \rar k$ of degree $n$ such that
$$ \langle v', w\rangle v'' \,=\, \pm \langle v, w''\rangle w'\,,\,\,\,\,\,\, \forall\,\,v,w\,\in\,C,$$
where $v'$ and $v''$ are the components of the coproduct of $v$ written in the Sweedler notation. A DG coalgebra $C$ is $n$-cyclic if it is $n$-cyclic as a graded coalgebra and
$$ \langle du, v \rangle \pm \langle u, dv \rangle \,=\, 0\,,$$
for all homogeneous $u,v \,\in\,{C}$, i.e, if $\langle \mbox{--}, \mbox{--} \rangle \,:\, {C}[n] \otimes C[n] \rar k[n]$ is a map of complexes. {\it By convention, we say that $C\,\in\,\DGC_{k/k}$ is $n$-cyclic if $\bar{C}$ is $n$-cyclic as a non-counital DG coalgebra}.

Assume that $C\,\in\,\DGC_{k/k}$ is equipped with a cyclic pairing of degree $n$ and let $R\,:=\,\cb(C)$ denote the (associative) cobar construction of $C$. Recall that $R\,\cong\, T_k(\bar{C}[-1])$  as a graded $k$-algebra. For $v_1,\ldots,v_n\,\in\,\bar{C}[-1]$, let $(v_1, \ldots ,v_n)$ denote the element $v_1 \otimes \ldots \otimes v_n$ of $R$. By~\cite[Theorem~15]{BCER}, the cyclic pairing on $C$ of degree $n$ induces a double Poisson bracket of degree $n+2$ (in the sense of~\cite{VdB})
$$\{\!\{ \mbox{--}, \mbox{--}\}\!\}\,:\, \bar{R} \otimes \bar{R} \rar {R} \otimes {R}$$
given by the formula
\begin{align} \la{dpbr}
\begin{aligned}
 &\{\!\{(v_1, \ldots, v_n), (w_1, \ldots, w_m)\}\!\}\,=\, \\
 &\sum_{\stackrel{i=1,\ldots,n}{j = 1,\ldots, m}} \pm \langle v_i, w_j \rangle  (w_1 ,\ldots, w_{j-1}, v_{i+1} ,\ldots, v_n)\otimes (v_1, \ldots, v_{i-1}, w_{j+1} ,\ldots, w_m) \,\text{.} \end{aligned} \end{align}
The above double bracket can be extended to $R \otimes R$ by setting $\{\!\{r,1\}\!\}\,=\, \{\!\{1,r\}\!\}=0$. Let $\{ \mbox{--}, \mbox{--}\}$ be the bracket associated to~\eqref{dpbr}:
\begin{equation} \la{bronr} \{ \mbox{--}, \mbox{--}\} \,:=\, \mu\,\circ\, \{\!\{ \mbox{--},\mbox{--}\}\!\}\,:\, {R} \otimes {R} \rar {R}\,,\end{equation}
 where $\mu$ is the multiplication map on ${R}$.  Let $\n\,:\, {R} \rar R_\n$ be the canonical projection and let $\{ \mbox{--}, \mbox{--}\}\,:\, \n \circ \{ \mbox{--}, \mbox{--}\}\,:\, {R} \otimes {R} \rar R_\n$. We recall that the bimodule ${R} \otimes {R}$ (with outer $R$-bimodule structure) has a double bracket (in the sense of~\cite[Defn. 3.5]{CEEY}) given by the formula
\begin{align*}
&\{\!\{ \mbox{--}, \mbox{--}\}\!\}\,\,:\,{R} \times ({R} \otimes {R}) \rar {R} \otimes ({R} \otimes {R}) \oplus ({R} \otimes {R}) \otimes {R}\,, \\
& \{\!\{r, p \otimes q\}\!\} \,:=\,  \{\!\{r,p\}\!\} \otimes q \oplus (-1)^{|p|(|r|+n)} p \otimes \{\!\{ r,q\}\!\} \,\text{.}
\end{align*}
This double bracket restricts to a double bracket on the sub-bimodule $\Omega^1R$ of $R \otimes R$ (~\cite[Corollary~5.2]{CEEY}). Let $\{\mbox{--}, \mbox{--}\}\,:\, R \otimes \Omega^1R \rar \Omega^1$ be the map $\mu \circ \{\!\{ \mbox{--},\mbox{--}\}\!\}$, where $\mu$ is the bimodule action map and let $\{\mbox{--}, \mbox{--}\}\, :\, R \otimes \Omega^1R \rar \Omega^1R_{\n}$ denote the map $\n \circ \{ \mbox{--},\mbox{--}\}$.

 The bracket $\{ \mbox{--},\mbox{--}\}\,:\, {R} \otimes {R} \rar R_\n$ descends to a DG $(n+2)$-Poisson structure on $R$. In particular, it descends to a (DG) Lie bracket $\{\mbox{--},\mbox{--}\}_{\n}$  on $R_\n$ of degree $n+2$. The restriction of the bracket~\eqref{bronr} to $\bar{R}$ induces a degree $n+2$ DG Lie module structure over $R_\n$ on $\bar{R}$ and the bracket $\{\mbox{--},\mbox{--}\}\,:\,R \otimes \Omega^1R \rar \Omega^1R_\n$ induces a degree $n+2$ DG Lie module structure over $R_\n$ on $\Omega^1R_\n$ (see~\cite[Proposition~3.11]{CEEY}). On homologies, we have (see~
\cite{CEEY}, Theorem~1.1 and Theorem~1.2)
\bthm \la{liestronhom}
Let $A\,\in\,\DGA_{k/k}$ be an augmented associative algebra Koszul dual to $C\,\in\,\DGC_{k/k}$. Assume that $C$ is $n$-cyclic. Then, \\
$(i)$ $\rHC_{\bullet}(A)$ has the structure of a graded Lie algebra (with Lie bracket of degree $n+2$).\\
$(ii)$ $\overline{\HH}_{\bullet}(A)$ has a graded Lie module structure over $\rHC_{\bullet}(A)$ of degree $n+2$.\\
$(iii)$ The maps $S,B$ and $I$ in the Connes periodicity sequence~\eqref{connessbi} are homomorphisms of degree $n+2$ graded Lie modules over $\rHC_{\bullet}(A)$.
\ethm
The Lie bracket of degree $n+2$ on $\rHC_{\bullet}(A)$ that is induced by a $(n+2)$-Poisson structure on $R_{\n}$ as above is an example of  a derived $(n+2)$-Poisson structure on $A$.\\

\noindent
\textbf{Convention.} Since we work with algebras that are Koszul dual to $n$-cyclic coalgebras, all Lie algebras that we work with have Lie bracket of degree $n+2$. Similarly, all Lie modules are degree $n+2$ Lie modules. We therefore, drop the prefix ``degree $n+2$" in the sections that follow. Following this convention, we shall refer to (derived) $(n+2)$-Poisson structures as (derived) Poisson structures.

\subsection{The Hodge filtration} \la{s2.2}
Consider the filtration on (the graded vector space) $\rHC_{\bullet}(\mathcal U\mfa)$ given by
\begin{equation}
\la{eqfilthc}
F_p\rHC_{\bullet}(\mathcal U\mfa)\,:=\, \bigoplus_{r \leq p+2} \HC^{(r)}_{\bullet}(\mfa) \,\text{.}\end{equation}
 Let $\mfa\,\in\,\DGL_k$ be Koszul dual to $C\,\in\,\cDGC_{k/k}$, where $C$ is $n$-cyclic. The following theorem is one of our main results.
\bthm \la{thodgefiltration}
The derived Poisson bracket $\{ \mbox{--},\mbox{--}\}$ on $\rHC_{\bullet}(\mathcal U\mfa)$ respects the filtration~\eqref{eqfilthc}.
Moreover,\\
\noindent
$(i)$ $\{ \HC^{(1)}_{\bullet}(\mfa), \HC^{(p)}_{\bullet}(\mfa)\} \,\subseteq \, \HC^{(p-1)}_{\bullet}(\mfa)$ for all $p \geq 1$. In particular, $\{\HC^{(1)}_{\bullet}(\mfa), \HC^{(1)}_{\bullet}(\mfa)\}=0$.\\
$(ii)$ $\{\HC^{(2)}_{\bullet}(\mfa), \HC^{(p)}_{\bullet}(\mfa)\} \subseteq \, \HC^{(p)}_{\bullet}(\mfa)$ for all $p \geq 1$.\\
In particular,  $\HC^{(2)}_{\bullet}(\mfa)$ is a Lie subalgebra of $\rHC_{\bullet}(\mathcal U\mfa)$ and $\rHC_{\bullet}(\mathcal U\mfa)$ is a Hodge weight graded Lie module over $\HC^{(2)}_{\bullet}(\mfa)$.
\ethm
Recall that $\mathcal L\,:=\, \cb_{\mathtt{Comm}}(C)$ gives a cofibrant resolution $\mathcal L \stackrel{\sim}{\rar} \mfa$  of $\mfa$ in $\DGL_k$ and the (associative) cobar construction $R\,:=\,\cb(C)$ gives a cofibrant resolution $R \stackrel{\sim}{\rar} \mathcal U\mfa$ of $\mathcal U\mfa$ (with $R\,\cong\, \mathcal U\mathcal L$). By Theorem~\ref{liestronhom}, the $n$-cyclic pairing on $C$ induces a derived Poisson structure on $\rHC_{\bullet}(\mathcal U\mfa)$. The proof of Theorem~\ref{thodgefiltration} is based on the following proposition. For notational brevity, let $V\,:=\, \bar{C}[-1]$.
\bprop \la{hodgebracketr}
Let $R^{(p)}\,:=\, \Sym^p(\mathcal L)$ as in Section~\ref{s1.3}. Let $\{ \mbox{--}, \mbox{--}\}$ be defined as in~\eqref{bronr}. \\
$(i)$  $\{ R^{(1)}, R^{(p)}\} \,\subseteq \, R^{(p-1)}$ for any $p \geq 1$. In particular, $\{R^{(1)}, R^{(1)}\}=0$.\\
$(ii)$ $ \{R^{(2)}, R^{(p)}\}\,\subseteq \,  R^{(p)}$ for any $p \geq 1$.\\
$(iii)$ For $p,q>2$, $\{R^{(q)}, R^{(p)}\} \,\subseteq\, \bigoplus_{r \leq p+q-2} R^{(r)}$.
\eprop
\bproof
Let $\alpha\,\in\,R^{(q)}$ be homogeneous. By~\cite[Section 2.4]{VdB}, the map $\{\alpha , \mbox{--}\}\,:\, \bar{R} \rar R$ is a derivation of degree $|\alpha|+n+2$. By Lemma~\ref{derivations}, it suffices to show that  $\{\alpha, V\}\,\subseteq \, \Sym^{q-1}(\mathcal L)$. Note that if $w\,\in\, V$, then by~\eqref{dpbr},
\begin{equation*} \{(v_1, \ldots, v_n), w\} \,=\, \sum_{i=1}^n \pm \langle v_i, w\rangle (v_{i+1} , \ldots, v_n, v_1, \ldots ,v_{i-1}) \,\text{.}\end{equation*}
It is easy to see that this coincides with the expression $\langle \bar{\partial}(v_1, \ldots, v_n), w \rangle$, where
$\langle \mbox{--}, \mbox{--} \rangle\,:\, \bar{R} \otimes V \rar R$ denotes the composite map
$$ \begin{diagram} \bar{R} \otimes V & \rTo^{\bar{\partial} \otimes \id_V} & (R \otimes V) \otimes V & \rTo & R \otimes (V \otimes V) & \rTo^{\id_V \otimes \langle \mbox{--}, \mbox{--} \rangle} & R \end{diagram} \,\text{.}$$
The unlabelled arrow in the above diagram swaps factors. Now, if $\alpha\,\in, R^{(q)}$, $\bar{\partial}(\alpha)\,\in\, \Sym^{q-1}(\mathcal L) \otimes V$ by Lemma~\ref{cycderhamhodge}. Thus, $\{\alpha, V\}\,\subseteq \, \Sym^{q-1}(\mathcal L)$ as desired.
\eproof
\bcor \la{hodgencpoiss} Let $\{\mbox{--},\mbox{--}\}_{\n}$ be the DG Poisson bracket on $R_{\n}$ as in Section~\ref{s2.1} above. Then,\\
$(i)$  $\{ R^{(1)}_\n, R^{(p)}_\n\}_{\n} \,\subseteq \, R^{(p-1)}_{\n}$ for any $p \geq 1$. In particular, $\{R^{(1)}_\n, R^{(1)}_\n\}_\n=0$.\\
$(ii)$ $ \{R^{(2)}_\n, R^{(p)}_{\n}\}_\n\,\subseteq \,  R^{(p)}_\n$ for any $p \geq 1$.\\
$(iii)$ For $p,q>2$, $\{R^{(q)}_\n, R^{(p)}_\n\}_\n \,\subseteq\, \bigoplus_{r \leq p+q-2} R^{(r)}_\n$.
\ecor
Corollary~\ref{hodgencpoiss} follows immediately from Proposition~\ref{hodgebracketr}. Theorem~\ref{thodgefiltration} follows immediately
from \eqref{hodgeds} and Corollary~\ref{hodgencpoiss}.

\subsection{Hodge filtration on Poisson modules} Let $\mfa,C, R, V$ and $\mathcal L$ be as in Section~\ref{s2.2}, with $C$ having a cyclic pairing of degree $n$. Recall that by Theorem~\ref{liestronhom}, $\overline{\HH}_{\bullet}(\mathcal U\mfa)$ is a Lie module over $\rHC_{\bullet}(\mathcal U\mfa)$.  In addition to~\eqref{eqfilthc}, we define the filtration on $\overline{\HH}_{\bullet}(\mathcal U\mfa)$
\begin{equation} \la{eqfilthh} F_p\overline{\HH}_{\bullet}(\mathcal U\mfa)\,:=\, \bigoplus_{r \leq p+2} \HH^{(r)}_{\bullet}(\mfa) \,\text{.} \end{equation}
and let $\{\mbox{--}, \mbox{--}\}\,:\, \rHC_{\bullet}(\mathcal U\mfa) \times \overline{\HH}_{\bullet}(\mathcal U\mfa) \rar \overline{\HH}_{\bullet}(\mathcal U\mfa)$ denote the action map.
\bthm \la{pstrhoch1}
 $\overline{\HH}_{\bullet}(\mathcal U\mfa)$ is a filtered Lie module, with filtration given by~\eqref{eqfilthh}. Moreover, for all $p \geq 0$,\\
$(i)$ $ \{\HC^{(1)}_{\bullet}(\mfa) , \HH^{(p)}_{\bullet}(\mfa)\} \subseteq \HH^{(p-1)}_{\bullet}(\mfa)$. \\
$(ii)$ $\{\HC^{(2)}_{\bullet}(\mfa),  \HH^{(p)}_{\bullet}(\mfa)\} \subseteq \HH^{(p)}_{\bullet}(\mfa)$.\\
In particular, $\overline{\HH}_{\bullet}(\mathcal U\mfa) $ is a Hodge weight graded Lie module over $\HC^{(2)}_{\bullet}(\mfa)$.
\ethm
 The restrictions of the maps $S,B$ and $I$ of the Connes periodicity sequence to the Hodge summand~\eqref{conneshodgep} give maps of graded vector spaces
\begin{align}
 \la{sp1} & S\,:\, \HC^{(p)}_{\bullet}(\mfa) \rar \HC^{(p+1)}_{\bullet-2}(\mfa)\\
 \la{bp1} & B\,:\, \HC^{(p+1)}_{\bullet}(\mfa) \rar \HH^{(p)}_{\bullet+1}(\mfa)\\
\la{ip1}   & I\,:\, \HH^{(p)}_{\bullet}(\mfa) \rar \HC^{(p)}_{\bullet}(\mfa) \,\text{.}
\end{align}
If we equip $\rHC_{\bullet}(\mathcal U\mfa)$ with filtration~\eqref{eqfilthc} and $\overline{\HH}_{\bullet}(\mathcal U\mfa)$   with filtration~\eqref{eqfilthh}, the maps $S,B$ and $I$ become filtered maps; more precisely, we have:
\begin{align}
\la{sp2} &  S\,:\, F_{\bullet}\rHC_{\bullet}(\mathcal U\mfa) \rar F_{\bullet+1}\rHC_{\bullet-2}(\mathcal U\mfa)\\
\la{bp2} &  B \,:\, F_{\bullet+1}\overline{\HC}_{\bullet}(\mathcal U\mfa) \rar F_{\bullet}\overline{\HH}_{\bullet+1}(\mathcal U\mfa)\\
\la{ip2} &  I\,:\, F_{\bullet}\overline{\HH}_{\bullet}(\mathcal U\mfa) \rar F_{\bullet}\overline{\HC}_{\bullet}(\mathcal U\mfa)\,\text{.}
\end{align}
The first statement in the following theorem is a refinement of~\cite[Theorem~1.2]{CEEY}.
\bthm \la{pstrhoch2}
With definitions~\eqref{sp2}-\eqref{ip2}, the maps $S,B$ and $I$ become {\rm filtered} Lie module maps. Moreover, in the Hodge summand~\eqref{conneshodgep}, these maps are module maps over the Lie algebra $\HC^{(2)}_{\bullet}(\mathcal U\mfa)$.
\ethm
\bcor \la{pstrGK}
The space $\mathrm{HA}_{\bullet}(\mathtt{Lie}, \mfa)$ has the structure of a graded Lie algebra, and $\mathrm{HB}_{\bullet}(\mathtt{Lie},\mfa)$ and $\mathrm{HC}_{\bullet}(\mathtt{Lie},\mfa)$ are graded Lie modules over $\mathrm{HA}_{\bullet}(\mathtt{Lie},\mfa)$.  Further, the maps $S,B$ and $I$ in the Connes periodicity sequence for (operadic) Lie cyclic homology (see Theorem~\ref{connesvsgksbi}) are maps of Lie modules.
\ecor
\bproof
Indeed, by definition, $\mathrm{HA}_{\bullet}(\mathtt{Lie},\mfa)\,=\, \HC^{(2)}_{\bullet}(\mfa)$. The remaining statements follow from Theorem~\ref{connesvsgksbi} and by putting $p=1$ in Theorem~\ref{thodgefiltration}, Theorem~\ref{pstrhoch1} and Theorem~\ref{pstrhoch2}.
\eproof
\subsubsection{Proof of Theorem~\ref{pstrhoch1}}

 Recall from Section~\ref{s1.3} that $\Omega^1R_\n\,\cong\, R \otimes V$ and that there is a direct sum decomposition (of complexes) $\Omega^1R_\n \,\cong\, \bigoplus_{p \geq 0}  \theta^{(p)}(\mathcal L)$. Further recall that  the isomorphism $\Omega^1R_\n\,\cong\, R \otimes V$ identifies $\theta^{(p)}(\mathcal L)$ with $R^{(p)} \otimes V$, where $R^{(p)}\,:=\, \Sym^p(\mathcal L)$. Recall from Section~\ref{s2.1} that there is a bracket $\{\mbox{--},\mbox{--}\}  \,:\, \bar{R} \times \Omega^1R_{\n} \rar \Omega^1R_{\n}$ inducing the structure of a DG Lie module over $R_\n$ on $\Omega^1R_{\n}$.
\bprop \la{hodgencpoissonforms} For any $p \geq 0$, the following inclusions hold:\\
$(i)$ $\{R^{(1)}, R^{(p)} \otimes V\} \subseteq R^{(p-1)} \otimes V$. In particular, $\{R^{(1)}, V\} =0$.\\
$(ii)$ $\{R^{(2)}, R^{(p)} \otimes V\} \subseteq R^{(p)} \otimes V$. Thus, $\Omega^1R_\n$ is a Hodge weight graded DG Lie module over $R^{(2)}_{\n}$.\\
$(iii)$ For any $m > 2$, $\{R^{(m)}, R^{(p)} \otimes V\} \subseteq \bigoplus_{ r \leq p+m-2} R^{(r)} \otimes V$.
  \eprop
\bproof
Let $r\,\in\,R^{(m)}$, $q\,\in\, R^{(p)}$ and $v \,\in \,V$ be homogeneous. By~\cite[Lemma 5.5]{CEEY},
\begin{equation} \la{ebforms} \{r, q \otimes v\}_{\Omega^1R_\n} \,=\, \n \circ [ \{r,q\}\cdot v \otimes 1 - \{r,q\} \otimes v + (-1)^{(|r|+n)|q|} (q\cdot \{r,v\} \otimes 1 -q \otimes \{r,v\})]\,,  \end{equation}
where $\Omega^1R$ on the right hand side is realized as a sub-bimodule of $R \otimes R$ equipped with the outer bimodule structure and where $\n\,:\,\Omega^1R \rar \Omega^1R_{\n}$ is the canonical projection. The first summand $\n \circ ( \{r,q\}\cdot v \otimes 1 - \{r,q\} \otimes v)$ of~\eqref{ebforms} is equal to $\{r,q\} \otimes v$ (after the identification of $\Omega^1R_\n$ with $R \otimes V$. The remaining (i.e, second) summand of~\eqref{ebforms} may be written as
$\n \circ [ q \cdot d\{r,v\}]$, where $d\,:\,R \rar \Omega^1R$ is the universal derivation.

If $m=1$, Proposition~\ref{hodgebracketr} implies that $\{r,q\}\,\in\, R^{(p-1)}$ and that $\{r,v\}\,\in\, k$. Thus, in this case, the second sumand vanishes while the first summand is in $R^{(p-1)} \otimes V$. This proves (i).

If $m=2$,  Proposition~\ref{hodgebracketr} implies that $\{r,q\}\,\in\, R^{(p)}$ and that $\{r,v\}\,\in\, \mathcal L$. Thus, the first statement of (ii) will follow once we show that for any $\alpha\,\in\,\mathcal L$, $\n \circ [ q \cdot d\alpha] \,\in\, R^{(p)} \otimes V$. Note that if $\n \circ [ q \cdot d\alpha] \,\in\, R^{(p)} \otimes V$ for all $q\,\in\,R^{(p)}$ for a given homogeneous $\alpha \in \mathcal L$, then for any $v \,\in\,V$ homogeneous,
\begin{align*}
\n \circ ( q \cdot d[v,\alpha])\,&=\, \n \circ (q \cdot d[v\cdot \alpha -(-1)^{|v||\alpha|} \alpha \cdot v])
\,=\, \n \circ (q \cdot [dv \cdot \alpha +v \cdot d\alpha -(-1)^{|v||\alpha|} (d\alpha \cdot v+ \alpha \cdot dv)])\\
\,&=\, \n \circ ((-1)^{|\alpha|(|q|+|v|)} [\alpha \cdot q \cdot dv - (-1)^{|q||\alpha|} q \cdot \alpha \cdot dv] + [q \cdot v -(-1)^{|q||v|} v \cdot q] \cdot d\alpha) \\
\,&=\, \n \circ ((-1)^{|\alpha|(|q|+|v|)} [\alpha, q] \cdot dv + [q,v] \cdot d\alpha)\,\text{.}
\end{align*}
Since $[\alpha,q]$ and $[q,v]$ are in $R^{(p)}$,  we have shown that $\n \circ ( q \cdot d[v,\alpha])$ is in $R^{(p)} \otimes V$ for all $q$. Thus, by induction on $n$, we show that $\n \circ ( q \cdot d[v,\alpha_n])$ is in $R^{(p)} \otimes V$ for $\alpha_n\,=\, [v_1,[v_2, \ldots [v_{n-1}, v_n] \ldots ]]$ for homogeneous $v_1,\ldots,v_n$ in $V$. This proves that $\n \circ [ q \cdot d\alpha] \,\in\, R^{(p)} \otimes V$ for all $q\,\in\,R^{(p)}$ for all $\alpha \in \mathcal L$ as desired. The second statement in (ii) follows from the first and the fact that the bracket on $\Omega^1R_\n$ descends to an DG Lie module structure on $\Omega^1R_\n$ over $R_\n$ (of which $R^{(2)}_{\n}$ is a DG Lie subalgebra).

Finally, for $m>2$, the first summand of~\eqref{ebforms}, which can be identified with $\{r,q\} \otimes v$, lies in $\bigoplus_{r \leq p+m-2} R^{(r)} \otimes V$, by Proposition~\ref{hodgebracketr}. Also note that $\{r,v\}\,\in\, R^{(m-1)}$, by Proposition~\ref{hodgebracketr}. Thus, the second summand $\n \circ [q \cdot d\{r,v\}]$ is of the form $\n \circ [q \cdot d\beta]$, where $\beta\,\in\,R^{(m-1)}$. This last expression is a $k$-linear combination of expressions of the form
$\n \circ [q \cdot \gamma \cdot d\alpha \cdot \phi] \,=\, \pm \n \circ  [\phi \cdot q \cdot \gamma \cdot d\alpha]$, where $\alpha\,\in\,\mathcal L$, and $\gamma\,\in\,R^{(i)}$ and $\phi\,\in\,R^{(j)}$ with $i+j=m-2$. Now, for any such expression, $\phi \cdot q \cdot \gamma\,\in\,  \bigoplus_{ r \leq p+m-2} R^{(r)}$.
That $ \n \circ  [\phi \cdot q \cdot \gamma \cdot d\alpha]$ lies in $ \bigoplus_{ r \leq p+m-2} R^{(r)} \otimes V$ now follows from the computation we made while proving (ii).
\eproof
It follows from Proposition~\ref{hodgebracketr}, Proposition~\ref{hodgencpoissonforms} and~\cite[Theorem~5.3]{CEEY} that the statement of Proposition~\ref{hodgencpoissonforms} holds word for word with $R^{(p)} \otimes V$ replaced by $\mathrm{Tot}\,X_2^{+,(p)}$. Theorem~\ref{pstrhoch1} then follows on homologies.
\subsubsection{Proof of Theorem~\ref{pstrhoch2}}
Let $F_pX^+(R)$ be given by $\oplus_{ r \leq p+2} X^{+, (r)}(\mfa)$. By Proposition~\ref{hodgebracketr}, Proposition~\ref{hodgencpoissonforms} and~\cite[Theorem~5.3]{CEEY}, there is an exact sequence of filtered DG Lie modules over $R_\n$
$$\begin{diagram} 0 & \rTo &  F_{\bullet}\mathrm{Tot}\,X_2^{+}(R) & \rTo^I & F_{\bullet}\mathrm{Tot}\,X^{+}(R)& \rTo^S & F_{\bullet+1}\mathrm{Tot}\,X^{+}(R)[2] & \rTo & 0 \end{diagram}\,\text{.}$$
This implies (i) on homologies. Again by Proposition~\ref{hodgebracketr}, Proposition~\ref{hodgencpoissonforms} and~\cite[Theorem~5.3]{CEEY}, the short exact sequence of complexes~\eqref{xcomplexhptot} is an exact sequence of DG Lie modules over $R^{(2)}_n$. This gives (ii).

\subsection{Examples}
\la{Sect3.5}
We illustrate the results of this section on several examples.

\subsubsection{Abelian Lie algebras} \la{s2.2.2} Let $\mfa=V$ be an abelian Lie algebra in homological degree $0$. Then, $A\,:=\,\mathcal U\mfa\,=\, \Sym(V)$. The Koszul dual cocommutative coalgebra is $C\,:=\,\Sym^c(V[1])$. By~\cite[Proposition~9.4]{CEEY}, $C$ has (upto scalars) a unique cyclic pairing $\langle \mbox{--},\mbox{--} \rangle$ that is of degree $-n$, where $n\,:=\,\dim_kV$. For a basis $\{v_1,\ldots,v_n\}$ of $V$, this pairing is uniquely determined by $\langle 1,v_1 \wedge \ldots \wedge v_n \rangle$.

In this case, there is a Hochschild-Kostant-Rosenberg (HKR) isomorphism $\rHC_{n}(A)\,\cong\, \Omega^n(A)/d\Omega^{n-1}(A)$. Further, there is an isomorphism of $A$-modules $\Psi\,:\, \Theta_p(A) \cong \Omega^{n-p}(A)$  given by $\Theta_p(\xi) \,=\, \iota_{\xi} \omega$, where $\omega$ is a fixed constant nonzero $n$-form and where $\iota_{\xi}$ is contraction by $\xi$.  In this case, by~\cite[Corollary 9.6]{CEEY}, the derived Poisson structure (of degree $2-n$) on $\rHC_{\bullet}(A)$ has a very explicit description after identifying reduced cyclic homology with forms via the HKR as above:
\begin{equation} \la{bracketofdiffforms} \{\alpha, \beta\} \,=\, (-1)^{(n-|\alpha|-1)(n-|\beta|)} \iota_{\eta}d\alpha\,,\,\,\,\, \text{ where } \eta\,:=\, \Psi^{-1}(d\beta)\,\text{.}\end{equation}
Now, by~\cite[Section 2]{bfprw2}, $\HC^{(p)}_{\bullet}(\mfa)$ is identified (via the HKR map) with $\Omega^{n,(p)}(A)/d\Omega^{n-1,(p+1)}(A)$, where $\Omega^{n,(p)}(A)$ is the space of $n$-forms on $A$ whose polynomial coefficients are homogeneous of weight $p$, i,e, $\Sym^p(V) \otimes \wedge^n(V)$. An easy computation using~\eqref{bracketofdiffforms} shows that $\{\HC^{(p)}_{\bullet}(\mfa), \HC^{(q)}_{\bullet}(\mfa)\}$ is actually contained in $\HC^{(p+q-2)}_{\bullet}(\mfa)$. Thus, in this case the derived Poisson Lie algebra from Theorem~\ref{thodgefiltration} becomes a graded (not just filtered) Lie algebra with respect to the Hodge decomposition.

\subsubsection{Necklace Lie algebras} Our next example is when $\mfa\,=\,LV$, the free Lie algebra generated by an even dimensional vector space (concentrated in homological degre $0$) equipped with a symplectic form $\langle \mbox{--}, \mbox{--} \rangle\,:\, V \times V \rar k$. Taking $\bar{C}:= V[1]$ with $0$ coproduct, we see that $C$ is Koszul dual to $\mfa$. Further, the symplectic form on $V$ can be viewed as a cyclic pairing of degree $-2$ on $\bar{C}$. The resulting derived Poisson bracket $\{\mbox{--},\mbox{--}\}$ equips $R_\n\,:=\, T_kV_\n$ with the structure of a Lie algebra. The Lie algebra $(R_{\n},  \{\mbox{--},\mbox{--}\})$ is isomorphic to the necklace Lie algebra (see~\cite{BL,G}) associated with a one vertex quiver having $\frac{1}{2}\dim_k V$ loops.

Theorem~\ref{thodgefiltration} implies that the necklace Lie bracket restricts to give a Lie bracket on the direct summand $\HC^{(2)}_0(\mfa)\,=\,\Sym^{(2)}(LV)_{\n}$ of $T_kV$.  Further, the filtration defined in Theorem~\ref{thodgefiltration} makes the necklace Lie algebra a filtered Lie algebra. This example also shows that the filtered structure in Theorem~\ref{thodgefiltration} is the best that we can get in general: for example, in the necklace Lie algebra generated by a two dimensional symplectic vector space $V$ with $v,w \,\in\,V$ such that $\langle v,w\rangle =1$, we have
$$\{ [v^3]_{\n}, [w^3]_{\n}\} \,=\, 9\,[v^2w^2]_{\n}\,\text{.}  $$
It is not difficult to verify that while the right hand side above is contained in $\oplus_{r \leq 4} \Sym^r(LV)_\n$, it is not contained in $\Sym^4(LV)_\n$.

Further, let $V_n\,=\,k^{2n}$ with the standard symplectic form. The symplectic form preserving inclusion $V_n \hookrightarrow V_{n+1}$ induces a  homomorphism $L_n \rar L_{n+1}$ of necklace Lie algebras, where $L_n\,:=\, ([TV_n]_\n, \{\mbox{--},\mbox{--}\})$. It is easily seen that these maps are compatible with the filtration from Theorem~\ref{thodgefiltration}. Thus, $L_{\infty}\,:=\, \varinjlim_n L_n$ admits a filtration as in Theorem~\ref{thodgefiltration}. In particular,  $L^{(2)}_{\infty}\,:=\, \varinjlim_n \Sym^2(LV_n)_\n$ is a Lie subalgebra of $L_{\infty}$.

Let $\mathcal P$ be a cyclic Koszul operad and let $W$ be a symplectic vector space with symplectic form $\omega$. Let $F_{\mathcal P}W$ be the free $\mathcal P$-algebra generated by $W$. Let  $\mathrm{Der}(F_{\mathcal P}W, \omega)$ be the Lie algebra of $\mathcal P$-algebra derivations of $F_{\mathcal P}W$ that preserve $\omega$ (see~\cite[Section 6]{G} for example). We refer to $\mathrm{Der}(F_{\mathcal P}W, \omega)$ as the Lie algebra associated with (the symplectic vector space) $W$ for the operad $\mathcal P$. It is not difficult to see that $\Sym^2(LV_n)_\n$ is the Lie algebra associated with $V_n$ for the Lie operad while $L_n$ is the Lie algebra associated with $V_n$ for the associative operad.  The Lie homology of $L_{\infty}$ is related to the cohomology of coarse moduli spaces of algebraic curves of fixed genus and fixed number of punctures while the homology of $L^{(2)}_{\infty}$ is related to spaces of outer automorphisms of free groups with punctures (see~\cite{Ko}).

\subsubsection{Unimodular Lie algebras} Let $\mfa$ be a Lie algebra of finite dimension $n$ in homological degree $0$. The Chevalley-Eilenberg coalgebra $C\,:=\,\C(\mfa;k)$ is Koszul dual to $\mfa$ and isomorphic to $\Sym^c(\mfa[1])$ as graded vector spaces. For the pairing $\langle \mbox{--},\mbox{--} \rangle$ from Section~\ref{s2.2.2} to be compatible with differentials, it is necessary and sufficient that $\mfa$ be {\it unimodular} i.e, that $\Tr(\mathrm{ad}(x))=0$ for all $x\,\in\,\mfa$. This is the case, for instance, if $\mfa$ is semisimple. Thus, for unimodular $\mfa$, $\rHC_{\bullet}(\mfa)$ has a derived Poisson structure of degree $2-n$, to which Theorem~\ref{thodgefiltration} applies. Given that the graded $(2-n)$-Lie structure on $\cb(C)_{\n}$ is identical to that of Section~\ref{s2.2.2} (with the Chevalley-Eilenberg diferential being the only new ingredient), it is reasonable to expect that for any unimodular Lie algebra $\mfa$, $\{\HC^{(p)}_{\bullet}(\mfa), \HC^{(q)}_{\bullet}(\mfa)\} \subseteq \HC^{(p+q-2)}_{\bullet}(\mfa)$ for all $p,q$.

\section{Topological applications}
\label{Sect4}

In this section, we will give a topological interpretation of Lie Hodge decompositions.
This interpretation is based on Quillen's famous theorem \cite{Q2} assigning to each simply connected topological space $ X $ a DG Lie algebra $\mfa_X$ called a Lie model of $X$. The cyclic homology of  $\, \U \mfa_X $ can be identified with the $S^1$-equivariant homology of the free loop space $ \LL X$ of $ X $, and our main observation is that the Hodge
components $ \HC^{(p)}_{\bullet}(\mathfrak{a}_X) $ correspond precisely to eigenspaces
of Frobenius maps  under this identification.
For an application, we will look at the Chas-Sullivan Lie algebra of a simply connected
closed manifold and show that the corresponding Lie bracket respects the Hodge filtration, thus making the string topology algebra a {\it filtered} Lie algebra.

Throughout this section, we assume that $k= \Q $ and all homology and cohomology groups are taken with rational coefficients.

\subsection{Hodge decomposition of $S^1$-equivariant homology}
\la{Sect4.1}
Let $X$ be a 1-connected 
topological space of finite rational type.
Recall ({\it cf.} \cite{FHT}) that one can associate to $ X $ a commutative
cochain DG algebra $ {\mathcal A}_X $, called a {\it Sullivan model} of $X$,
and a connected (chain) DG Lie algebra $ \mfa_X $, called a {\it Quillen model} of $X$.
Each of these algebras is uniquely determined up to homotopy and each encodes the
rational homotopy type of $X$. The relation between
them is given by a DG algebra quasi-isomorphism
\begin{equation}
\la{SQ}
\C^{\bullet}(\mathfrak{a}_X;\mathbb Q) \stackrel{\sim}{\to} {\mathcal A}_X\ ,
\end{equation}
where $\, \C^{\bullet}(\mathfrak{a}_X;\Q) \,$ is the Chevalley-Eilenberg cochain complex of $ \mfa_X $.

Now, let $\LL X$ denote the free loop space of $X$, i.e. the space of all continuous
maps $ S^1 \to X $ equipped with compact open topology. This space carries a natural
$ S^1$-action (induced by rotations of $ S^1 $), and  one  can define the $S^1$-equivariant
homology of $ \LL X $ via the Borel construction:
$$
\H^{S^1}_{\bullet}(\LL X) := \H_\bullet(ES^1 \times_{S^1} \LL X, \Q)\ .
$$
We will use a  {\it reduced} version of equivariant homology, which is defined by
$$
\rH^{S^1}_{\bullet}(\LL X) := \Ker[\,\H^{S^1}_{\bullet}(\LL X) \xrightarrow{\pi_*} \H_\bullet(BS^1)\,]\ ,
$$
where the map $ \pi_* $ comes from the natural homotopy fibration
\begin{equation}
\la{fibr}
\LL X \to ES^1 \times_{S^1} \LL X \xrightarrow{\pi} BS^1 \ .
\end{equation}
The following theorem is a well known result due to Goodwillie \cite{Go} and Jones \cite{J}
(see also \cite{JM}).
\bthm[\cite{J}]
\la{top1}
There are natural isomorphisms of graded vector spaces
\begin{equation*}
\la{ax}
\alpha_X:\, {\rHH}_{\bullet}(\mathcal U\mathfrak{a}_X)  \xrightarrow{\sim} \overline{\H}_{\bullet}(\mathcal L X)
\ ,\qquad
\beta_X:\, \rHC_{\bullet}(\mathcal U\mathfrak{a}_X) \xrightarrow{\sim} {\rH}^{S^1}_{\bullet}(\LL X)
\end{equation*}
transforming the Connes periodicity sequence for $ \, \U\mathfrak{a}_X $ to the Gysin long exact sequence for the $S^1$-equivariant homology of $ \LL X $.
\ethm

Next, for each integer $n \ge 0 $, we consider the $n$-fold covering of the circle:
$$
\omega^n\,:\,S^1 \rar S^1\ ,\quad e^{i\theta} \mapsto e^{in\theta}\ ,
$$
and denote by $\,\varphi_X^n\,;\,\mathcal L X \rar \mathcal L X\,$ the induced map on $ \LL X$. While the maps $\varphi_X^n$ are not equivariant with respect to the $S^1$-action on $ \LL X$, they fit into the commutative diagram ({\it cf.} \cite{BFG})
\begin{equation}\la{diag}
\begin{diagram}
\LL X & \rTo^{\varphi_X^n} & \LL X\\
 \dTo & & \dTo\\
(ES^1 \times_{S^1} \mathcal L X)_\Q & \rTo^{\tilde{\varphi}_X^n} &  (ES^1 \times_{S^1} \mathcal L X)_\Q \\
  \dTo^{} & & \dTo_{}\\
  (BS^1)_\Q & \rTo^{B\omega^n} & (BS^1)_\Q\\
  \end{diagram}\ ,
\end{equation}
where the columns arise from the rationalization (Bousfield localization at zero) of the homotopy fibration \eqref{fibr}. The maps $ \varphi_X^n $ and $ \tilde{\varphi}_X^n $ in \eqref{diag}  induce graded linear endomorphisms
\begin{equation*}
\Phi_X^n:\,  {\rH}_{\bullet}(\LL X) \to  {\rH}_{\bullet}(\LL X)
\ ,\qquad \tilde{\Phi}_X^n :\,  {\rH}^{S^1}_{\bullet}(\LL X)  \to   {\rH}^{S^1}_{\bullet}(\LL X)\ ,
\end{equation*}
which are called the {\it power} or {\it Frobenius operations} on $ {\rH}_{\bullet}(\LL X) $  and $ {\rH}^{S^1}_{\bullet}(\LL X) $. We write
\begin{equation}\la{phi}
{\rH}^{(d)}_{\bullet}(\LL X)\, := \bigcap_{n \ge 0}\,\Ker(\Phi_X^n - n^d\,\id)\ ,\qquad
{\rH}^{S^1,\, (d)}_{\bullet}(\LL X)\, :=
 \bigcap_{n \ge 0}\,\Ker(\tilde{\Phi}_X^n - n^d\,\id)
\end{equation}
for the common eigenspaces of these endomorphisms corresponding to the eigenvalues $ n^d $.

Now, the main result of this section can be stated as follows.
\bthm
\la{top2}
For each $p\ge 0$, there are natural  isomorphisms
\begin{equation*}
\overline{\HH}^{(p)}_{\bullet}(\mathfrak{a}_X) \cong {\rH}^{(p)}_{\bullet}(\mathcal L X) \ ,\qquad
\rHC^{(p)}_{\bullet}(\mathfrak{a}_X)  \cong
 {\rH}^{S^1, \,(p-1)}_{\bullet}(\LL X)
\end{equation*}
given by restriction of the isomorphisms $ \alpha_X $ and $\beta_X$  of Theorem~\ref{top1}.
\ethm
To prove Theorem~\ref{top2} we recall from \cite{bfprw1} the construction of a Lie Hodge decompostion in terms of (dual) Adams operations on the cobar construction of a cocommutative DG algebra. First, recall that by Adams
operations on a commutative algebra $A$ one usually means a
family $\{\psi^n\}_{n\ge 0}$ of algebra homomorphisms
$ \psi^n: A \to A $ satisfying the relations
\begin{equation}
\la{adams1}
\psi^1 = \id\ ,\quad \psi^n \circ \psi^m = \psi^{nm}\ .
\end{equation}
Dually, for a cocommutative coalgebra $ A $,
we define Adams operations to be a
family $\{\psi^n\}_{n\ge 0}$ of {\it coalgebra} homomorphisms
$ \psi^n:  A \to   A $ satisfying \eqref{adams1}.
Now, if $ A $ is either a commutative or  cocommutative Hopf algebra, there are natural Adams operations on $A$ defined by the formulas\footnote{For commutative Hopf algebras, these Adams operations are defined and studied in \cite[Section~4.5]{L}, while, for cocommutative Hopf algebras, they appear in \cite[Section~7]{bfprw1}.}
\begin{equation}
\la{adams2}
\psi^1 = \id\ ,\quad \psi^n = \mu_n \circ \Delta^n\ ,\ n \ge 2\ ,
\end{equation}
where $\,\mu_n:\,A^{\otimes n} \to A\,$ is the $n$-fold product and $\,\Delta^n: A \to  A^{\otimes n}\, $ is the  $n$-fold coproduct
on $A$.

Next, recall that, if $ \mfa $ is a DG Lie algebra, its Koszul dual is  a cocommutative (conilpotent) DG coalgebra $C$. Associated to $C$
are two cobar constructions: the classical (Adams) cobar construction
$ \bOmega(C) $, which gives a cofibrant resolution of $ \U\mfa $ in $ \DGA_{k/k} $, and the Lie cobar construction $ \bOmega_{\tt Comm}(C) $, which gives a cofibrant resolution of $ \mfa $  in $ \DGL_k $. The relation between these constructions is given by a canonical isomorphism $ \bOmega(C) \cong \U[ \bOmega_{\tt Comm}(C)] $, which shows that $ \bOmega(C) $ has a natural structure of a cocommutative DG Hopf algebra. Thus, $ \bOmega(C) $ can be equipped with a  family $ \{\psi^n\}_{n\ge 0} $ of (coalgebra) Adams operations given by formulas \eqref{adams2}. It is shown in \cite{bfprw1} (see {\it op. cit.}, Section~7) that these operations on $ \bOmega(C) $ induce Adams operations on the (reduced) cyclic homology of $\, \U\mfa $, and the construction of \cite{bfprw1} can be extended to define Adams operations on the Hochschild homology of $\, \U\mfa $.
We denote these Adams operations on $\, \rHH_\bullet(\U \mfa) \,$ and
$\, \rHC_\bullet(\U \mfa) \,$ by $ \Psi^n $ and $ \tilde{\Psi}^n $, respectively.
Proposition~7.3 of \cite{bfprw1} then implies that the Hodge components $ \HH_\bullet^{(p)}(\mfa) $ and $ \HC_\bullet^{(p)}(\mfa) $ are  the common
eigenspaces of $ \Psi^n $ and $ \tilde{\Psi}^n $ corresponding to the eigenvalues $ n^p $:
\begin{equation}
\la{psi}
{\rHH}_\bullet^{(p)}(\mfa)\, =
 \bigcap_{n \ge 0}\,\Ker(\Psi^n - n^p\,\id)\ , \qquad
{\rHC}_\bullet^{(p)}(\mfa)\, =
 \bigcap_{n \ge 0}\,\Ker(\tilde{\Psi}^n - n^p\,\id)\ .
\end{equation}

Now, let $ \Psi_X^n $ and $ \tilde{\Psi}_X^n $ denote the Adams operations on $\, \rHH_\bullet(\U \mfa_X) $ and $\, \rHC_\bullet(\U\mfa_X) $ coming from  $\cb(C) $,
where $ C = \C_{\bullet}(\mathfrak{a}_X; \Q) $ is the Chevalley-Eilenberg chain complex of $ \mfa_X $. Theorem~\ref{top2} follows immediately from  \eqref{phi}, \eqref{psi} and the next key proposition.
\bprop
\la{intert}
For each $ n \ge 0 $,  the isomorphism $\alpha_X $ intertwines the Adams operation $ \Psi_X^n $ with the Frobenius operation $ \Phi_X^n $, and the isomorphism $ \beta_X$ intertwines $ \tilde{\Psi}_X^n $ with
the $n$-th multiple of $ \tilde{\Phi}^n_X$, i.e.
$$
\alpha_X \circ \Psi^n_X \,=\, \Phi_X^n \circ \alpha_X \ ,\qquad
\beta_X \circ \tilde{\Psi}^n_X \, = \, n\,\tilde{\Phi}^n_X \circ \beta_X\ .
$$
\eprop
\bproof
We deduce Proposition~\ref{intert} from results of  the paper \cite{BFG}. Theorem~B of that paper provides natural isomorphisms
\begin{equation} \la{abx}
a_X:\, {\rHH}_{- \bullet}({\mathcal A}_X)  \xrightarrow{\sim} {\rH}^{\bullet}(\mathcal L X)
\ ,\qquad
b_X:\, \rHC_{- \bullet}({\mathcal A}_X) \xrightarrow{\sim} {\rH}_{S^1}^{\bullet - 1}(\LL X)\ ,
\end{equation}
relating the Hochschild and cyclic homology of a Sullivan model of
$ X $ to the (reduced) cohomology and $ S^1 $-equivariant  cohomology of $ \LL X $. Since  $ {\mathcal A}_X $ is a commutative DG algebra,
its Hochschild and cyclic homology carry natural Adams operations
which we denote by $ \Psi_n $ and $ \tilde{\Psi}_n $, respectively. On the other hand, the cohomology and $ S^1 $-equivariant  cohomology of $ \LL X $ carry Frobenius operations induced by $ \varphi^n $: we denote these by $ \Phi_n:\,  {\rH}^{\bullet}(\mathcal L X) \to  {\rH}^{\bullet}(\mathcal L X)  $ and $ \tilde{\Phi}_n:
{\rH}_{S^1}^{\bullet}(\LL X) \to {\rH}_{S^1}^{\bullet}(\LL X) $.
Now, part $(2)$ of \cite[Theorem~B]{BFG} says that
\begin{equation}
\la{ab3}
a_X \circ \Psi_n \,=\, \frac{1}{n}\, \Phi_n \circ a_X \ ,\qquad
b_X \circ \tilde{\Psi}_n \, = \, \tilde{\Phi}_n \circ b_X\ .
\end{equation}
To relate the isomorphisms \eqref{abx} to those of Theorem~\ref{top1}, we recall a theorem of Quillen~\cite{Q} identifying
\begin{equation} \la{quil}
{\rHH}_{\bullet}(\mathcal U\mathfrak{a}_X)\,\cong\,  \overline{\HH}_{\bullet}[\C_\bullet(\mathfrak{a}_X;\mathbb Q)] \ ,
\qquad \rHC_{\bullet}(\mathcal U\mathfrak{a}_X) \,\cong\, \rHC_{\bullet+1}[\C_\bullet(\mathfrak{a}_X;\mathbb Q)]\ .
\end{equation}
Since $X$ is assumed to be of finite rational type, its Lie and Sullivan models $ \mfa_X $ and $ {\mathcal A}_X $ are locally finite DG algebras, i.e. have finite-dimensional components in each homological degree. Dualizing
\eqref{SQ}, we then have a quasi-isomorphism of DG coalgebras $\, {\mathcal A}_X^* \stackrel{\sim}{\to} \C^{\bullet}(\mathfrak{a}_X;\Q)^* \cong \C_{\bullet}(\mathfrak{a}_X;\Q)\,$.
Combined with \eqref{abx} and \eqref{quil}, this quasi-isomorphism induces natural isomorphisms
\begin{eqnarray}
&& \overline{\HH}_{\bullet}(\mathcal U\mathfrak{a}_X) \,\cong\,  \overline{\HH}_{\bullet}[\C_\bullet(\mathfrak{a}_X;\mathbb Q)] \,\cong\,  {\rHH}_{-\bullet}({\mathcal A}_X)^{\ast} \,\xrightarrow{(a_X^{-1})^*}\, \overline{\H}^{\bullet}(\mathcal L X)^{\ast} \,\cong\,  \overline{\H}_{\bullet}(\mathcal L X)\ , \la{mm1}\\
&&\rHC_{\bullet}(\mathcal U\mathfrak{a}_X) \,\cong\, \rHC_{\bullet+1}[\C_\bullet (\mathfrak{a}_X;\mathbb Q)] \,\cong\, \rHC_{-\bullet-1}({\mathcal A}_X)^{\ast} \,\xrightarrow{(b_X^{-1})^*}\,  \overline{\H}^{\bullet}_{S^1}(\mathcal L_X)^{\ast} \,\cong\, \overline{\H}^{S^1}_{\bullet}(\mathcal L X)\ ,\la{mm2}
\end{eqnarray}
where the star $(\,\mbox{--}\,)^{\ast}$ stands for graded linear duals.  Since $a_X$ and $b_X$ are functorial and transform Connes'  exact sequence on homology to the Gysin sequence for $S^1$-equivariant cohomology, the isomorphisms \eqref{mm1} and \eqref{mm2} are functorial and transform Connes' sequence to  the Gysin sequence on $S^1$-equivariant homology. Thus \eqref{mm1} and \eqref{mm2} coincide with the isomorphisms $\alpha_X$ and $\beta_X$ of Theorem~\ref{top1}.
Now, by~\cite[Proposition 7.4]{bfprw1}, the direct summand $\HC_{\bullet}^{(p)}(\mathfrak{a}_X)$ of $\rHC_{\bullet}(\mathcal U\mathfrak{a}_X)$ is identified with $ \rHC_{-\bullet-1}^{(p-1)}({\mathcal A}_X)^{\ast}$ under the first two isomorphisms in \eqref{mm2}, and by \cite[Theorem~B]{BFG}, $\, \rHC_{-\bullet-1}^{(p-1)}({\mathcal A}_X)^{\ast}$ is the common eigenspace of the Adams operations $  \tilde{\Psi}^*_n $ with the eigenvalues $n^{p-1}$. On the other hand, $ \HC_{\bullet}^{(p)}(\mathfrak{a}_X) $ is the common eigenspace of the Adams operations $\tilde{\Psi}^n $ with the eigenvalues $ n^p $. Hence the first two isomorphisms of \eqref{mm2} intertwine $\tilde{\Psi}^n $ on $\rHC_{\bullet}(\mathcal U\mathfrak{a}_X)$  with $\, n \tilde{\Psi}^*_n \,$ on $\, \rHC_{-\bullet-1}^{(p-1)}({\mathcal A}_X)^{\ast}$. Now, by \eqref{ab3}, the map $\, (b_X^{-1})^{\ast}$ intertwines $ \tilde{\Psi}^*_n $ with $\tilde{\Phi}^*_n $, and under the last isomorphism in \eqref{mm2}, the dual power operation $\tilde{\Phi}^*_n $ on $ \rH^{\bullet}_{S^1}(\mathcal L_X)^{\ast} $ corresponds to $ \tilde{\Phi}^n $ on $ \rH^{S^1}_{\bullet}(\mathcal L X) $. Thus,
the composite map \eqref{mm2} intertwines $\tilde{\Psi}^n $ with $\,n\,\tilde{\Phi}^n $. This proves the second equality of Proposition~\ref{intert}. The proof of the first is similar.
\eproof
Theorem~\ref{top2} together with Theorem~\ref{hodgesbi} has the following important corollary.
\bcor \la{hsstring}
There are natural Hodge decompositions
$$
{\rH}_{\bullet}(\mathcal L X) \,\cong\, \bigoplus_{p=0}^{\infty} \overline{\H}^{(p)}_{\bullet}(\mathcal L X)\ ,\qquad
\rH^{S^1}_{\bullet}(\mathcal L X)\,\cong\,\bigoplus_{p=0}^{\infty} \overline{\H}^{S^1,\,(p)}_{\bullet}(\mathcal L X)\ .
$$
The Gysin sequence decomposes into a direct sum of Hodge components
$$
 \ldots\, \xrightarrow{D} \,\overline{\H}^{S^1,\, (p+1)}_{n-1}(\mathcal L X) \,\to\,\overline{\H}^{(p)}_{n}(\mathcal L X) \,\to\, \overline{\H}^{S^1, (p)}_{n}(\mathcal L X) \,\xrightarrow{D}\, \overline{\H}^{S^1, \, (p+1)}_{n-2}(\mathcal L X)\,\to\, \ldots
$$
where $D$ stands for the Gysin map.
\ecor

\subsection{Hodge filtration on string topology}
Let $M$ be a simply connected closed oriented manifold of dimension $d$.
A construction of Lambrechts and Stanley (see \cite{LS}, Theorem~1.1) provides a finite-dimensional  commutative DG algebra $ \mathcal A$, which is a model for the singular cochain complex of $M$. This model comes equipped with a nondegenerate cyclic pairing of (cohomological) degree $ n = -d $. The linear dual of $ \mathcal A$ gives a cocommutative coalgebra model
$ C := {\mathcal A}^* $ for the singular chain complex of $M$ that has a cyclic pairing of (homological) degree $ n = -d $. It is known (and easy to check) that
$\cb_{\mathtt{Comm}}(C)$ is a Lie model of $M$; hence, by Theorem~\ref{top1}, we have
$$
\rHC_{\bullet}[\cb(C)]\,\cong\, \overline{\H}^{S^1}_{\bullet}(\mathcal L M)\,,\,\,\,\,\, \overline{\HH}_{\bullet}[\cb(C)]\,\cong\, \overline{\H}_{\bullet}(\mathcal L M) \ ,
$$
and the cyclic pairing on $C$ induces a derived Poisson structure on $\overline{\H}^{S^1}_{\bullet}(\mathcal L M)$ of degree $n+2$. The corresponding Lie bracket coincides with the Chas-Sullivan bracket \cite{ChS}, and the resulting Lie algebra is called the {\it string topology Lie algebra} of $M$ ({\it cf.}~\cite{CEG}). Further,  by Theorem~\ref{liestronhom}, $\overline{\H}_{\bullet}(\mathcal L M)$ is a Lie module (of degree $n+2$) over the string topology Lie algebra of $M$, and the Gysin map $D\,:\, \overline{\H}^{S^1}_{\bullet}(\mathcal L M) \rar \overline{\H}^{S^1}_{\bullet-2}(\mathcal L M)$ is a Lie module homomorphism.

We may now apply in this situation the results of Section~\ref{Sect3} and Section~\ref{Sect4}:
as a consequence of Theorem~\ref{thodgefiltration}, Theorem~\ref{pstrhoch1}, Theorem~\ref{pstrhoch2} and Theorem~\ref{top2}, we get
\bthm \la{tstringhomology}
$(i)$
The string topology Lie algebra of a closed $d$-dimensional manifold $M$ is filtered as a Lie algebra with respect to the following {\rm Hodge filtration}
$$
F_p\overline{\H}^{S^1}_{\bullet}(\mathcal L M) \,:=\, \bigoplus_{q \leq p+1} \overline{\H}^{S^1,\, (q)}_{\bullet}(\mathcal L M)\,\text{.} $$

\noindent
$(ii)$ The homology of the free loop space $\overline{\H}_{\bullet}(\mathcal L M)$ is filtered as a Lie module over the string topology Lie algebra of $M$  with respect to the following Hodge filtration
$$ F_p\overline{\H}_{\bullet}(\mathcal L M) \,:=\, \bigoplus_{q \leq p+2} \overline{\H}^{(q)}_{\bullet}(\mathcal L M)\,\text{.} $$

\noindent
$(iii)$ The Chas-Sullivan bracket restricts to the first Hodge component $\, {\rH}^{S^1,\,(1)}_{\bullet}(\mathcal L M)$, making it a Lie algebra. Further, $ {\rH}^{S^1}_{\bullet}(\mathcal L M)$  is a graded Lie module over $\overline{\H}^{S^1,\,(1)}_{\bullet}(\mathcal L M)$ with the grading given by the Hodge decomposition of $\overline{\H}^{S^1}_{\bullet}(\mathcal L M)$.

\noindent
$(iv)$ For each $ p \ge 0 $, the Gysin map $\,D:\, \overline{\H}^{S^1,\,(p)}_{\bullet}(\mathcal L M) \rar \overline{\H}^{S^1,\,(p+1)}_{\bullet-2}(\mathcal L M)$  is a map of graded Lie modules over $\,\overline{\H}^{S^1, (1)}_{\bullet}(\mathcal L M)$.
\ethm
The Lambrechts-Stanley Lie model $\cb_{\mathtt{Comm}}({\mathcal A}^*)$ and unimodular Lie algebras (in particular, abelian Lie algebras) are Koszul dual to ``Poincar\'{e} duality CDGC'', i.e. cocommutative DG coalgebras that are graded linear duals of Poincar\'{e} duality CDGA's in the sense of~\cite{LS}. This should be contrasted with
necklace Lie algebras, in which case the coproduct on the Koszul dual coalgebra and the cyclic pairing are unrelated. Given the example of symmetric algebras in Section~\ref{s2.2.2}, where the Hodge filtration actually becomes a Hodge decomposition, we expect the following conjecture to be true.
\bconj
\la{conj1}
Let $\mfa\,\in\, \DGL_k$ be Koszul dual to a connected finite-dimensional $C\,\in\,\cDGC_{k/k}$ equipped with a Poincar\'{e} duality pairing of degree $n$. Then, the Hodge filtrations on the derived Poisson structures on $\rHC_{\bullet}(\mfa)$ and $\overline{\HH}_{\bullet}(\mfa)$ become direct Hodge decompositions, i.e. for all $p,q \ge 0 $,
$$\{\HC^{(p)}_{\bullet}(\mfa),\, \HC^{(q)}_{\bullet}(\mfa)\} \,\subseteq \, \HC^{(p+q-2)}_{\bullet}(\mfa)\ ,\quad \{\HC^{(p)}_{\bullet}(\mfa),\, \HH^{(q)}_{\bullet}(\mfa)\} \,\subseteq \, \HH^{(p+q-2)}_{\bullet}(\mfa)\ .
$$
In particular, we expect  that the Chas-Sullivan bracket of a closed $d$-dimensional manifold satisfies
$$
\{\overline{\H}^{S^1,\,(p)}_{\bullet}(\mathcal L M), \ \overline{\H}^{S^1,\,(q)}_{\bullet}(\mathcal L M)\}\, \subseteq\, \overline{\H}^{S^1,\, (p+q-1)}_{\bullet}(\mathcal L M)\ .
$$
\econj

\section{Relation to derived representation schemes} \la{s4}

In this section, we recall the construction of derived representation schemes of Lie algebras from~\cite[Section 6, Section 7]{bfprw1}. The main result of this section is Theorem~\ref{dreppoiss}, the proof of which we outline. The full details will appear in~\cite{Yin}.

\subsection{Derived representation schemes and Drinfeld traces}

 Let $\g$ be a finite dimensional Lie algebra. Consider the functor
$$(\mbox{--})_{\g} \,:\, \DGL_k \rar \cDGA_{k/k}\,, \,\,\,\,\mfa \,\mapsto \, \mfa_{\g}\,:=\, \frac{\Sym_k(\mfa \otimes \g^{\ast})}{\langle\langle ( x \otimes \xi_1).(y \otimes \xi_2) -(y \otimes \xi_1).(x \otimes \xi_2) -[x,y] \otimes \xi \rangle \rangle} \,,$$
where $\g^{\ast}$ is the vector space dual to $\g$ and where $\xi \mapsto \xi_1 \wedge \xi_2$ is the map dual to the Lie bracket on $\g$. The augmentation on $\mfa_\g$ is the one induced by the map taking the generators $\mfa \otimes \g^{\ast}$ to $0$. Let $\g({\mbox{--}})\,:\, \cDGA_{k/k} \rar \DGL_k$ denote the functor $B \mapsto \g(\bar{B})\,:=\, \g \otimes \bar{B}$. It is shown in~\cite[Section 6.3]{bfprw1} that the functors $(\mbox{--})_{\g}\,:\, \DGL_k \rightleftarrows \cDGA_{k/k}\,:\, \g({\mbox{--}})$ form a (Quillen) adjoint pair.

Thus, $\mfa_\g$ is the commutative (DG) algebra corresponding to the (DG) scheme $\Rep_\g(\mfa)$ parametrizing representations of $\mfa$ in $\g$. Since the functor $(\mbox{--})_{\g}$ is left Quillen, it has a well behaved left derived functor
$$\L(\mbox{--})_{\g}\,:\,\Ho(\DGL_k) \rar \Ho(\cDGA_{k/k})\,\text{.}$$
Like for any left derived functor, we have  $\L(\mfa)_\g\,\cong\, \mathcal L_\g $ in $\Ho(\cDGA_{k/k})$, where $\mathcal L \stackrel{\sim}{\rar} \mfa$ is any cofibrant resolution in $\DGL_k$. We define
$$\DRep_\g(\mfa)\,:=\, \L(\mfa)_\g \, \text{ in } \Ho(\cDGA_{k/k})\,,\,\,\,\,\H_{\bullet}(\mfa,\g)\,:=\, \H_{\bullet}[\L(\mfa)_\g]\,\text{.}$$
$\DRep_\g(\mfa)$ is called the {\it derived representation algebra} for representations of $\mfa$ in $\g$. The homology $\H_{\bullet}(\mfa,\g)$ is called the {\it representation homology} of $\mfa$ in $\g$. It is not difficult to check that $\g$ acts naturally by derivations on the graded (commutative) algebra $\H_{\bullet}(\mfa,\g)$. We denote the corresponding (graded) subalgebra of $\g$-invariants by $\H_{\bullet}(\mfa,\g)^{\ad\,\g}$.

Let $\mathcal L \stackrel{\sim}{\rar} \mfa$ be a cofibrant resolution. The unit of the adjunction $(\mbox{--})_{\g}\,:\, \DGL_k \rightleftarrows \cDGA_{k/k}\,:\, \g({\mbox{--}})$ is the universal representation
$$ \pi_\g\,:\, \mathcal L \rar \g(\mathcal L_\g) \,\text{.}$$
Let the functor $\lambda^{(p)}$ be as in Section~\ref{s1.1}.  There is a natural map $\lambda^{(p)}[\g(\mathcal L_\g)] \rar \mathcal L_\g \otimes \lambda^{(p)}(\g)$. For $P\,\in\, I^p(\g)\,:=\,\Sym^p(\g^{\ast})^{\ad\,\g}$, evaluation at $P$ gives a linear functional $\mathrm{ev}_P$ on $\lambda^{(p)}(\g)$. One thus has the composite map
$$\begin{diagram} \lambda^{(p)}(\mathcal L) & \rTo^{\lambda^{(p)}(\pi_\g)} &  \lambda^{(p)}[\g(\mathcal L_\g)]  & \rTo&  \mathcal L_\g \otimes \lambda^{(p)}(\g) & \rTo^{\id \otimes \mathrm{ev}_P} & \mathcal L_\g \end{diagram} $$
for $P\,\in\,I^p(\g)$. On homologies, this gives the map
$$ \Tr_\g(P, \mfa)\,:\, \HC^{(p)}_{\bullet}(\mfa) \rar \H_{\bullet}(\mfa,\g)^{\ad\,\g}\,,$$
which we call the {\it Drinfeld trace map } associated to $P$ (see~\cite[Section 7]{bfprw1} for further details regarding this construction). If $\g$ is semisimple, the Killing form is a canonical element of $I^2(\g)$. We denote the associated Drinfeld trace  by
$$\Tr_\g(\mfa)\,:\,\HC^{(2)}_{\bullet}(\mfa) \rar \H_{\bullet}(\mfa,\g)^{\ad\,\g}\,\text{.}$$
Let $\mfa\,\in\,\DGL_k$ be Koszul dual to $C\,\in\,\cDGC_{k/k}$, where $C$ is $n$-cyclic.
\bthm \la{dreppoiss}
 There is a Poisson structure on $\H_{\bullet}(\mfa,\g)^{\ad\,\g}$ such that the Drinfeld trace map $\Tr_\g(\mfa)$ is a graded Lie algebra homomorphism.
\ethm
Theorem~\ref{dreppoiss} may thus be viewed as a generalization of~\cite[Theorem 6.7]{G} for the Lie operad (the latter result is proven only for free algebras over operads). It could also be seen as an analog of~\cite[Theorem 2]{BCER} in the Lie setting.

\subsection{Proof of Theorem~\ref{dreppoiss}}

Recall from~\cite[Section 4.5]{GK} that a cyclic pairing $\langle \mbox{--}, \mbox{--} \rangle$ on a (graded) Lie algebra $\mfa$ is a symmetric pairing that is $\ad$-invariant. Equivalently, the map $\langle \mbox{--}, \mbox{--} \rangle \,:\, \mfa \otimes \mfa \rar k$ is a (graded) $\mfa$-module homomorphism, where $k$ is equipped with the trivial action. Dually, a cyclic pairing on a (graded) Lie coalgebra $\mathfrak{G}$ is a pairing $\langle \mbox{--}, \mbox{--} \rangle\,:\, \mathfrak{G} \otimes \mathfrak{G} \rar k$ that is a (graded) $\mathfrak{G}$-comodule homomorphism, with $k$ equipped with the trivial coaction. Explicitly, if $\langle \mbox{--},\mbox{--} \rangle$ is a cyclic pairing of degree $n$ on a Lie coalgebra $\mathfrak{G}$, then if $\Delta(x) \,=\, x^1 \otimes x^2$, etc. in the Sweedler notation, we have
$$x^1 \langle x^2, y \rangle \pm x^2 \langle x^1, y \rangle \pm y^1 \langle x, y^2 \rangle \pm y^2 \langle x, y^1 \rangle \,=\,0\,, $$
where the signs are determined by the Koszul sign rule. If $\mathfrak{G}$ is differential graded, we further demand that a cyclic pairing of degree $n$ on $\mathfrak{G}$ be compatible with differential, i.e, that for all $x,y\,\in\,\mathfrak{G}$, $\langle \delta x, y \rangle +(-1)^{|x|+n}\langle x, \delta y \rangle =0$.

Further recall that for a DG Lie coalgebra $\mathfrak{G}$, one has the Chevalley-Eilenberg algebra $\C^c(\mathfrak{G};k)$ which is the construction formally dual to the Chevalley-Eilenberg coalgebra $\C(\mfa;k)$ os a DG Lie algebra.
\blemma \la{poissce}
If $\mathfrak{G}\,\in\,\DGLC_k$ is equipped with a cyclic pairing of degree $n$, then the Chevalley-Eilenberg algebra $\C^c(\mathfrak{G};k)$ acquires a DG Poisson structure.
\elemma
\bproof
Note that as a graded commutative algebra, $\C^c(\mathfrak{G};k)\,\cong\, \Sym(\mathfrak{G}[-1])$. The degree $n$-cyclic pairing on $\mathfrak{G}$ gives a skew symmetric degree pairing on $\mathfrak{G}[-1]$. This gives a graded Poisson structure on $\Sym(\mathfrak{G}[-1])$. The compatibility of this structure with respect to the differential on  $\C^c(\mathfrak{G};k)$ follows from the fact that the pairing on $\mathfrak{G}$ is cyclic.
\eproof
Let $\mfa\,\in\,\DGL_k$ be Koszul dual to $C\,\in\,\cDGC_{k/k}$ where $C$ is $n$-cyclic. If $\g$ is semisimple, then the Killing form on $\g$ gives an isomorphism $\g \, \cong\,\g^{\ast}$. Under this isomorphism, the Killing form on $\g$ is identified with a cyclic pairing $\kappa$ on $\g^{\ast}$. Note that $\g^{\ast}(\bar{C})\,:=\,\g^{\ast} \otimes \bar{C}$ has the structure of a DG Lie coalgebra. Tensoring $\kappa$ with the pairing on $\bar{C}$, we obtain a degree $n$-cyclic pairing on $\g^{\ast}(\bar{C})$. By Lemma~\ref{poissce}, this gives a DG Poisson structure of degree on $\C^c(\g^{\ast}(C);k)$, which represents $\DRep_\g(\mfa)$ in $\Ho(\cDGA_{k/k})$ by~\cite[Theorem 6.5]{bfprw1}.  On homologies, we obtain a graded Poisson structure on $\H_{\bullet}(\mfa,\g)$. It can be shown that the above Poisson structure on $\H_{\bullet}(\mfa,\g)$  restricts to a Poisson structure on $\H_{\bullet}(\mfa,\g)^{\ad\,\g}$.

On the other hand,  the Killing form on $\g$ is an element of $I^2(\g)$. The associated Drinfeld trace is a map of graded vector spaces
$\Tr_\g(\mfa)\,:\, \HC^{(2)}_{\bullet}(\mfa) \rar \H_{\bullet}(\mfa,\g)^{\ad\,\g}$. Note that by Theorem~\ref{thodgefiltration} (ii), $\HC^{(2)}_{\bullet}(\mfa)$ acquires the structure of a graded Lie algebra, namely, its derived Poisson structure. On the other hand, by Lemma~\ref{poissce}, the pairing on $\bar{C}$ induces a graded Poisson structure on $\H_{\bullet}(\mfa,\g)$. We now study the relation between these Poisson structures.

\subsubsection{} Let $\mathcal L\,:=\, \cb_{\mathtt{Comm}}(C)$. Thus, by~\cite[Theorem 6.5]{bfprw1}, $\mathcal L_\g\,\cong\, \C^c(\g^{\ast}(\bar{C});k)$. By Proposition~\ref{hodgebracketr} (ii), $\lambda^{(2)}(\mathcal L)$ is a DG Lie algebra. In addition, by Lemma~\ref{poissce}, $\mathcal L_\g$ has a DG Poisson structure. Let $\Tr_{\g}(\mathcal L)$ be the Drinfeld trace $\lambda^{(2)}(\mathcal L) \rar \mathcal L_\g$ associated with the Killing form on $\g$. Let $R\,:=\,\cb(C)$, which is isomorphic to $\mathcal U\mathcal L$ as a DG algebra.  By Lemma~\ref{cycderhamhodge} and~\cite[Lemma A.1]{BR}, the map $\bar{\partial}\,:\, R \rar \Omega^1R_\n$ (see~\eqref{cyclicderham}) induces a map $\bar{\partial}\,:\, \lambda^{(2)}(\mathcal L) \rar \mathcal L \otimes V$, where $V\,:=\,\bar{C}[-1]$. Note that the DG-module of $\Omega^1_{\mathcal L_\g}$ of K\"{a}hler differentials on $\mathcal L_\g$ is isomorphic to $\mathcal L_\g \otimes \g^{\ast}(\bar{C})[-1]$ as a graded $\mathcal L_\g$-module.  Let $d\,:\,\mathcal L_\g \rar \Omega^1_{\mathcal L_\g}$ be the de Rham differential. The prof of the following lemma will appear in~\cite{Yin}.
\blemma \la{twooneforms}
The following diagram commutes:
$$
\begin{diagram}
\lambda^{(2)}(\mathcal L) & \rTo^{\bar{\partial}} & \mathcal L \otimes V & \rTo^{\pi_{\g} \otimes \id} & \mathcal L_\g \otimes \g \otimes V\\
  & \rdTo_{\Tr_\g(\mathcal L)} & & & \dTo_{\cong}  \\
& & \mathcal L_\g & \rTo^d & \Omega^1_{\mathcal L_\g}\\
     \end{diagram}
$$
Here, the vertical isomorphism on the right identifies $\g$ with $\g^{\ast}$ through the Killing form.
\elemma
It is not difficult to verify that the bracket on $\lambda^{(2)}(\mathcal L)$ is given by the composite map
$$ \begin{diagram} \lambda^{(2)}(\mathcal L) \otimes \lambda^{(2)}(\mathcal L) & \rTo^{\bar{\partial} \otimes \bar{\partial}} & (\mathcal L \otimes V) \otimes (\mathcal L \otimes V) & \rTo & (\mathcal L \otimes \mathcal L) \otimes (V \otimes V) & \rTo^{\mathrm{can} \otimes \langle \mbox{--}, \mbox{--} \rangle } & \lambda^{(2)}(\mathcal L) \end{diagram} \,,$$
where $\mathrm{can}\,:\, \mathcal L \otimes \mathcal L \rar \lambda^{(2)}(\mathcal L)$ is the canonical projection. By the construction of $\Tr_\g(\mathcal L)$, the following diagram commutes.
$$
\begin{diagram}
 (\mathcal L \otimes V) \otimes (\mathcal L \otimes V)  & \rTo^{(\pi_\g \otimes \id) \otimes (\pi_\g \otimes \id)} & (\mathcal L_\g \otimes \g(\bar{C})[-1] )  \otimes (\mathcal L_\g \otimes \g(\bar{C})[-1] )\\
\dTo & & \dTo\\
(\mathcal L \otimes \mathcal L) \otimes (V \otimes V) & & (\mathcal L_\g \otimes \mathcal L_\g) \otimes (\g^{\ast}(\bar{C})[-1] \otimes \g^{\ast}(\bar{C})[-1])\\
\dTo^{\mathrm{can} \otimes \langle \mbox{--}, \mbox{--} \rangle} & & \dTo_{\mu \otimes  \langle \mbox{--}, \mbox{--} \rangle}\\
 \lambda^{(2)}(\mathcal L_\g) & \rTo^{\Tr_\g} & \mathcal L_\g\\
\end{diagram}
 $$
Here, in the first vertical arrow on the right, $\g$ is identified with $\g^{\ast}$ via Killing form. $\mu$ is the product on $\mathcal L_\g$ and the pairing on $\g^{\ast}(\bar{C})[-1]$ is the one induced by the cyclic pairing on $\g^{\ast}(\bar{C})$. Composing the bottom arrow with the left vertical arrows composed with $\bar{\partial} \otimes \bar{\partial}$, we obtain $\Tr_\g(\{\mbox{--},\mbox{--}\})\,:\, \lambda^{(2)}(\mathcal L) \otimes \lambda^{(2)}(\mathcal L) \rar \mathcal L_\g$. On the other hand, by Lemma~\ref{twooneforms}  composing the bottom arrow with the left vertical arrows composed with $\bar{\partial} \otimes \bar{\partial}$ gives us the map
$$  \lambda^{(2)}(\mathcal L) \otimes \lambda^{(2)}(\mathcal L) \rar \mathcal L_\g\,,\,\,\,\, (\alpha,\beta) \mapsto \langle d\Tr_\g(\alpha), d\Tr_\g(\beta) \rangle \,,$$
where the pairing $\langle \mbox{--},\mbox{--} \rangle \,:\, \Omega^1_{\mathcal L_\g} \otimes \Omega^1_{\mathcal L_\g} \rar \mathcal L_\g$ is induced by the pairing on $\g^{\ast}(\bar{C})[-1]$. This map is easily seen to be $\{\Tr_\g(\mbox{--}),\Tr_\g(\mbox{--})\}\,:\,\lambda^{(2)}(\mathcal L) \otimes \lambda^{(2)}(\mathcal L) \rar \mathcal L_\g$. Thus,
$$\Tr_\g(\{\mbox{--},\mbox{--}\})\,=\, \{\Tr_\g(\mbox{--}),\Tr_\g(\mbox{--})\}\,:\,\lambda^{(2)}(\mathcal L) \otimes \lambda^{(2)}(\mathcal L) \rar \mathcal L_\g\,\text{.} $$
It is not difficult to verify that the DG Poisson bracket on $\mathcal L_\g$ restricts to a DG Poisson bracket on $\mathcal L_\g^{\ad\,\g}$. We thus obtain Theorem~\ref{dreppoiss} on homologies.

\section{Appendix: Free (Lie) algebras}
Let $V$ be a (homologically) graded $k$-vector space. Let $R\,:=\,T_kV$, the free graded $k$-algebra generated by $V$. Let $\mathcal L\,:=\, LV$, the free graded Lie algebra generated by $V$. It is well known that $R\,\cong\,\mathcal U\mathcal L$. Thus, the direct sum of symmetrization maps gives an isomorphism of graded $\mathcal L$-modules
$$ \Sym(\mathcal L) \,=\,\oplus_p \Sym^p(\mathcal L) \,\cong\, R \,\text{.}$$
In what follows, we shall view $\Sym^p(\mathcal L)$ as a graded subspace of $R$ via the symmetrization map. For $\beta_1,\ldots, \beta_p$ homogeneous in $R$, let $s(\beta_1,\ldots,\beta_p)$ denote the {\it symmetrization} of $\beta_1,\ldots, \beta_p$, i.e, the sum
$$s(\beta_1, \ldots,\beta_p) \,:=\, \sum_{ \sigma\,\in\,S_p} \pm \beta_{\sigma(1)} \cdot \ldots \cdot \beta_{\sigma(p)}\,,$$
where $\cdot$ denotes the product on $R$ and where the sign in front of each summand is the sign determined by the Koszul sign rule.
\blemma \la{derivations}
Let $\delta\,:\, R \rar R$ be a homogeneous derivation. Suppose that $\delta(V) \,\subseteq\, \Sym^q(\mathcal L)$.\\
(i) If $q>1$, then $\delta[\Sym^p(\mathcal L)] \,\subseteq\, \bigoplus_{r \leq p+q-1} \Sym^{r}(\mathcal L)$ for any $p \geq 1$.\\
(ii) If $q \leq 1$, then  $\delta[\Sym^p(\mathcal L)] \,\subseteq\, \Sym^{p+q-1}(\mathcal L)$.
\elemma
\bproof
First, note that if $\alpha\,\in\,\mathcal L$ is homogeneous and if $\delta(\alpha)\,\in\, \Sym^q(\mathcal L)$, then for any homogeneous $v\,\in\,V$,
\begin{align*}
\delta([\alpha, v])\,&=\, \delta(\alpha \cdot v -(-1)^{|\alpha||v|} v \cdot \alpha)\\
&\,=\, \delta(\alpha) \cdot v +(-1)^{|\alpha||\delta|} \alpha \cdot \delta(v) -(-1)^{|\alpha||v|} \delta(v) \cdot \alpha -(-1)^{|\alpha||v|+|\delta||v|} v \cdot \delta(\alpha)\\
&\,=\, [\delta(\alpha), v] +(-1)^{|\alpha||\delta|} [\alpha, \delta(v)] \,\text{.}\\
\end{align*}
Since $[\mathcal L, \Sym^q(\mathcal L)] \,=\, [\Sym^q(\mathcal L), \mathcal L] \,\subseteq \, \Sym^q(\mathcal L)$ and $V \subseteq \mathcal L$, the element $\alpha_n := [v_1,[v_2, \ldots [v_{n-1}, v_n] \ldots ]]$ satisfies $\delta(\alpha_n) \subseteq \Sym^q(\mathcal L)$ for homogeneous $v_1,\ldots,v_n$ in $V$ by induction on $n$. It follows that $\delta(\mathcal L)\,\subseteq\,\Sym^q(\mathcal L)$. The desired lemma now follows from the fact that for $\beta_1,\ldots,\beta_p\,\in\,\mathcal L$,
$$ \delta[s(\beta_1,\ldots,\beta_p)]\,=\, \sum_{i=1}^p \pm s(\beta_1,\ldots, \delta(\beta_i),\ldots,\beta_p) \,\text{.}$$
Indeed, the symmetrization of $p-1$ elements of $\mathcal L$ with an element of $\Sym^q(\mathcal L)$ is in $\bigoplus_{r \leq p+q-1} \Sym^{r}(\mathcal L)$ if $q>1$ and in $\Sym^{p+q-1}(\mathcal L)$ if $q=0,1$.
\eproof
Identify $\Omega^1R_\n$ with $R \otimes V$ as in~\eqref{omeganat}. Let $\bar{\partial}\,:\, \bar{R} \rar \Omega^1R_{\n}$ be as in~\eqref{cyclicderham}.
\blemma \la{cycderhamhodge}
$\bar{\partial} [\Sym^p(\mathcal L)] \,\subseteq\, \Sym^{p-1}(\mathcal L) \otimes V$.
\elemma
\bproof
Let $\bar{\partial}'\,:\, \bar{R} \rar V \otimes \bar{R}$ be the operator
$$ (v_1,\ldots, v_n) \mapsto \sum_{i=1}^n (-1)^{(|v_1|+\ldots+|v_{i-1}|)(|v_i|+\ldots+|v_n|)} v_i \otimes (v_{i+1} ,\ldots, v_n,v_1, \ldots ,v_{i-1}) \,\text{.}$$
Note that there is a {\it right} action of $S_n$ on $V^{\otimes n}$. Explicitly, for $\sigma\,\in\,S_n$ and for $v_1, \ldots, v_n$ homogeneous, $(v_1 \ldots v_n) \cdot \sigma \,:=\, \pm (v_{\sigma(1)} \ldots v_{\sigma(n)})$, where the sign is determined by the Koszul sign rule. Let $\tau$ be the $n$-cycle $(1 2 \ldots n)$.  Then, the restriction of $\bar{\partial}'$ to $V^{\otimes n}$ is given by the composite map
$$ \begin{diagram} V^{\otimes n} & \rTo^{(\mbox{--}) \cdot N} & V^{\otimes n} & \rTo & V \otimes V^{\otimes n-1} \end{diagram}\,,$$
where $N\,=\, \sum_{i=0}^{n-1} \tau^i$ and where the last arrow is the obvious isomorphism that permutes no factors. Further note that the above right action of $S_n$ on $V^{\otimes n}$ is dual to the left action of $S_n$ on $W^{\otimes n}$ used in~\cite[Section 4.5]{L} (with $W\,:=\,V^{\ast}$).

Let $S_{n,p}$ be the set of permutations in $S_n$ having $p-1$ descents in the sense of~\cite[Section 4.5.5]{L}. Let $l^P_n\,:=\, \sum_{\sigma\,\in\,S_{n,p}} \sigma$ and let $e^{(p)}_n$ be the Eulerian idempotent
$$e^{(p)}_n \,=\, \sum_{j=1}^n a^{p,j}_n l^j_n\,, $$
where the Stirling numbers $a^{p,j}_n$ are defined by the identity $\sum_{p=1}^n a^{p,j}_n X^p\,=\, \binom{X-j+n}{n}$.

By~\cite[Remark 2.10]{L2}, the right action of the Eulerian idempotent $e^{(p)}_n$ on $V^{\otimes n}$ is the projection from $V^{\otimes n}$ to $\Sym^p(\mathcal L) \cap V^{\otimes n}$. Let $S_{n-1}$ be viewed as the subgroup of $S_n$ fixing $1$. By the proof of~\cite[Theorem 4.6.6]{L} (more specifically, formula $(4.6.6.2)$ in {\it loc. cit.}),
$e^{(p)}_nN\,=\, Ne^{(p-1)}_{n-1}$. Hence, for any $\alpha\,\in\, \Sym^p(\mathcal L) \cap V^{\otimes n}$,
\begin{eqnarray*}
\alpha \cdot N \,=\, (\alpha \cdot e^{(p)}_n) \cdot N\,=\, \alpha \cdot (e^{(p)}_nN) \,=\, \alpha \cdot (Ne^{(p-1)}_{n-1}) \,=\, (\alpha \cdot N) \cdot e^{(p-1)}_{n-1}\,\text{.}
\end{eqnarray*}
It follows that $(\alpha \cdot N) \cdot e^{(i)}_{n-1}\,=\, \delta_{i,p-1} \alpha \cdot N$, where $\delta_{i,j}$ is the Kronecker delta. Thus,
$$\bar{\partial}' [\Sym^p(\mathcal L) \cap V^{\otimes n}] \,\subseteq\, V \otimes (\Sym^{p-1}(\mathcal L) \cap V^{\otimes n-1})\,\text{.} $$
The desired lemma follows once we observe that $\bar{\partial}$ is given by composing $\bar{\partial}'$ with the isomorphism $V \otimes R \,\cong\, {R} \otimes V$ that swaps factors.
\eproof

\subsection*{Acknowledgements}{\footnotesize
We would like to thank Mike Mandell for interesting discussions and suggestions.
Research of the first two authors was partially supported by the Simons Foundation Collaboration Grant 066274-00002B (`Representation Homology').}

\end{document}